%% file: main.tex
\documentclass[cis]{ipart}
\input{math_commands.tex}

\RequirePackage{hyperref}
\pdfoutput=1
\usepackage{graphicx}
\usepackage[utf8]{inputenc} 
\usepackage[T1]{fontenc}    
\usepackage{subfigure}
\usepackage{float}
\usepackage{natbib}
\usepackage{amsmath}
\usepackage{tabu}
\usepackage{booktabs}
\usepackage{siunitx}
\usepackage{multirow,multicol}
\usepackage{xcolor}

\usepackage{amsthm}

\usepackage{booktabs,dcolumn,caption}
\newcolumntype{d}[1]{D{.}{.}{#1}} 

\usepackage{graphicx}
\usepackage{graphics}
\usepackage{paralist}
\numberwithin{equation}{subsection}
\usepackage{wrapfig}

\newcommand\lizhen[1]{\textcolor{red}}

\newcommand*{\scale}[2][4]{\scalebox{#1}{$#2$}}

\usepackage{tikz}
\usetikzlibrary{matrix,decorations.pathreplacing}

\pubyear{2021}
\volume{0}
\issue{0}
\firstpage{1}
\lastpage{9}
\arxiv{2009.03892}

\begin{document}

\title[Neural-PDE: A RNN based neural network for solving time dependent PDEs ]{Neural-PDE: A RNN based neural network for solving time dependent PDEs }

\begin{aug}
    \author{\fnms{Yihao} \snm{Hu}\ead[label=e1]{yhu5@nd.edu}},
    \address{Department of Applied and Computational Mathematics and Statistics \\
         University of Notre Dame,
         Notre Dame, IN 46545, USA \\
             \printead{e1}}
    \author{\fnms{Tong} \snm{Zhao}\ead[label=e2]{tzhao2@nd.edu}},
    \address{Department of Computer Science and Engineering \\
         University of Notre Dame,
         Notre Dame, IN 46545, USA \\
             \printead{e2}}
    \author{\fnms{Shixin} \snm{Xu}\ead[label=e3]{shixin.xu@dukekunshan.edu.cn}},
    \address{Duke Kunshan University \\
Kunshan, Jiangsu 215316, P.R. China \\
             \printead{e3}}
    \author{\fnms{Lizhen} \snm{Lin}\ead[label=e4]{lizhen.lin@nd.edu}},
    \address{Department of Applied and Computational Mathematics and Statistics\\
         University of Notre Dame,
         Notre Dame, IN 46545, USA \\
             \printead{e4}}
    \and 
    \author{\fnms{Zhiliang} \snm{Xu}
            \ead[label=e5]{zxu2@nd.edu}}
    \address{Department of Applied and Computational Mathematics and Statistics\\
         University of Notre Dame,
         Notre Dame, IN 46545, USA \\
             \printead{e5}}
\end{aug}

\begin{abstract}
Partial differential equations (PDEs) play a crucial role in studying a vast number of problems in science and engineering. Numerically solving nonlinear and/or high-dimensional PDEs is frequently a challenging task. 
Inspired by the traditional finite difference and finite elements methods and emerging advancements in machine learning, we propose a  sequence-to-sequence learning (Seq2Seq) framework called Neural-PDE, which allows one to automatically learn governing rules of any time-dependent PDE system from existing data by using a bidirectional LSTM encoder, and predict the solutions in next $n$ time steps. One critical feature of our proposed framework is that the Neural-PDE is able to  simultaneously learn and simulate all variables of interest in a PDE system. We test the Neural-PDE by a range of  examples, from one-dimensional PDEs to a multi-dimensional and nonlinear complex fluids model. The results show that the Neural-PDE is capable of learning the initial conditions, boundary conditions and differential operators defining the initial-boundary-value problem of a PDE system without the knowledge of the specific form of the PDE system. In our experiments, the Neural-PDE can efficiently extract the dynamics within $20$ epochs training and produce accurate predictions.
Furthermore, unlike the traditional machine learning approaches for learning PDEs, such as CNN and MLP, which require great quantity of parameters for model precision, the Neural-PDE shares parameters among all time steps, and thus considerably reduces computational complexity and leads to a fast learning algorithm.
\end{abstract}
\maketitle


\section{Introduction}
\label{sec:introduction}
The research of time-dependent partial differential equations (PDEs) is regarded as one of the most important disciplines in applied mathematics. PDEs appear ubiquitously  in a broad spectrum of fields including physics, biology, chemistry,  and finance, to name a few. 
Despite their fundamental importance,  most PDEs  can not be solved analytically and have to rely on numerical solving methods. 
Developing efficient and accurate numerical schemes for solving PDEs, therefore, has been  an active research area over the past few decades~\citep{courant1967partial,osher1988fronts,leveque1992numerical,cockburn2012discontinuous,thomas2013numerical,johnson2012numerical}.   Still, devising  stable and accurate schemes with acceptable computational cost is a difficult task, especially when nonlinear and(or) high-dimensional PDEs are considered. Additionally, 
PDE models emerged from science and engineering disciplines usually require huge empirical data for model calibration and validation, and determining the multi-dimensional parameters in such a PDE system poses another challenge~\citep{peng2020multiscale}.

Deep learning is considered to be the state-of-the-art tool in classification and prediction of nonlinear inputs, such as image, text, and speech~\citep{litjens2017survey, devlin2018bert,lecun1998gradient,krizhevsky2012imagenet,hinton2012deep}.  Recently, considerable efforts have been made to  employ deep learning tools in designing data-driven methods for solving PDEs~\citep{Han2018SolvingHP,long2018pde,sirignano2018dgm,raissi2019physics}. Most of these approaches are based on fully-connected neural networks (FCNNs), convolutional neural networks(CNNs) and multilayer perceptron (MLP). These neural network structures usually require an increment of the layers to improve the predictive accuracy~\citep{raissi2019physics}, and subsequently lead to a more complicated model due to the additional parameters.     
Recurrent neural networks (RNNs) are another type of neural network architectures. RNNs predict the next time step value by using the input data from the current and previous states and share parameters across all inputs. 
This idea~\citep{sherstinsky2020fundamentals} of using current and previous step states to calculate the state at the next time step is not unique to RNNs. In fact, it is ubiquitously used in numerical PDEs. Almost all time-stepping numerical methods applied to solve time-dependent PDEs, such as Euler's, Crank-Nicolson, high-order Taylor and its variance Runge-Kutta~\citep{ascher1997implicit} time-stepping methods, update numerical solution by utilizing solution from previous steps. 
  

    

This motivates us to think what would happen if we  replace the previous step data in the neural network with numerical solution data to PDE supported on grids. It is possible that the neural network behaves like a time-stepping method, for example, forward Euler's method yielding the numerical solution at a new time point as the current state output~\citep{chen2018neural}. Since the  numerical solution on each of the grid point (for finite difference) or grid cell (for finite element) computed at a set of contiguous time points can be treated as neural network input in the form of one time sequence of data, the deep learning framework can be trained to predict any time-dependent PDEs from the time series data supported on some grids if the bidirectional structure is applied~\citep{huang2015bidirectional,schuster1997bidirectional}. In other words, the supervised training process can be regarded as a practice of the deep learning framework to learn the numerical solution from the input data, by learning the coefficients on neural network layers.

Long Short-Term Memory (LSTM)~\citep{hochreiter1997long} is a neural network built upon RNNs. Unlike vanilla RNNs, which suffer from losing long term information and high probability of gradient vanishing or exploding, LSTM has a specifically designed memory cell with a set of new gates such as input gate and forget gate. Equipped with these new gates which control the time to preserve and pass the information, LSTM is capable of learning long term dependencies without the danger of having gradient vanishing or exploding. In the past two decades, LSTM has been widely used in the field of natural language processing (NLP), such as machine translation, dialogue systems, question answering systems~\citep{lipton2015critical}.


Inspired by numerical PDE schemes and LSTM neural network, we propose a new deep learning framework, denoted as Neural-PDE. It simulates multi-dimensional governing laws, represented by time-dependent PDEs, from time series data generated on some grids and predicts the next $n$ time steps data. The Neural-PDE is capable of intelligently processing related data from all spatial  grids by using the bidirectional~\citep{schuster1997bidirectional} neural network, and thus guarantees the accuracy of the numerical solution and the feasibility in learning any time-dependent PDEs. The detailed structures of the Neural-PDE and data normalization are introduced in Section~\ref{sec:method}.

The rest of the paper is organized as follows. 
Section~\ref{sec:prelim} briefly reviews finite difference method and finite element method for solving PDEs.
Section~\ref{sec:method} contains detailed description of designing the Neural-PDE. 
In Section~\ref{sec:experiments}, we apply the Neural-PDE to solve four different PDEs, including the $1$-dimensional(1D) wave equation, the $2$-dimensional(2D) heat equation, and two systems of PDEs: the invicid Burgers' equations and a coupled Navier Stokes-Cahn Hilliard equations, which widely appear in multiscale modeling of complex fluid systems. We demonstrate the robustness of the Neural-PDE, which achieves accuracy within 20 epochs with an admissible mean squared error, even when we add Gaussian noise in the input data.

\section{Preliminaries}
\label{sec:prelim}
\subsection{Time Dependent Partial Differential Equations}
\label{sec:FD_FEM}
A time-dependent partial differential equation is an equation of the form:
\begin{equation}
    u_{t} = f(x_{1},\cdots, u, \frac{\partial u}{\partial x_{1}},\cdots,\frac{\partial u}{\partial x_{n}},\frac{\partial ^2 u}{\partial x_{1} \partial x_{1}},\cdots,\frac{\partial ^2 u}{\partial x_{1}\partial x_{n}}, \cdots,
    \frac{\partial ^n u}{\partial x_{1} \cdots \partial x_{n}})~,
\label{eq:time_pde_model}    
\end{equation}
where $u = u(x_1,...,x_n, t)$ is known,  $x_{i}\in \mathbb{R} $ are spatial variables, and the operator $f$ maps $\mathbb{R}^N \mapsto \mathbb{R}$. For example, consider the  parabolic heat equation: $u_{t} = \alpha^2 \Delta u$, where $u$ represents the temperature and $f$ is the Laplacian operator $\Delta$. 
Eq.~(\ref{eq:time_pde_model}) can be solved by  finite difference  methods, which are briefly reviewed below for the self-completeness of the paper.

\subsection{Finite Difference Method}
\label{sec:fdm}

Consider using a finite difference method (FDM) to solve a two-dimensional second-order PDE of the form:
\begin{equation}
    u_{t} = f(x, y, u_{x}, u_{y},u_{xx},u_{yy}), \quad (x,y) \in \Omega \subset \mathbb{R}^2, \quad t \in \mathbb{R}^{+} \cup\{0\}~,
\label{eq:FD_model_PDE}
\end{equation}
with some proper boundary conditions.
Let
$\Omega$ be $\Omega = [x_a,x_b]\times[y_a, y_b]$, and 
\begin{align}
    u_{i,j}^{n} = u( x_{i},y_{j}, t_{n})
\end{align}
where $t_{n} = n  \delta t, \; 0\leq n \leq N$, and  $\delta t = \frac{T}{N} \;$ for some large integer $N$. $x_{i} = i  \delta x, \; 0\leq i \leq N_x$, $\delta x = \frac{x_a - x_b}{N_x } \;$.~ $y_{j} = j  \delta y, \; 0\leq j \leq N_y$, $\delta y = \frac{y_a - y_b}{N_y } \;$. $N_x$ and $N_y$ are integers.

The central difference methods approximate the spatial derivatives as follows~\citep{thomas2013numerical}:
\begin{align}
    u_{x}(x_i, y_j, t) &= \frac{1}{2 \delta x}(u_{i+1,j} -u_{i-1,j}) + \mathcal{O}(\delta x^2)~,\\
     u_{y}(x_i, y_j, t) &= \frac{1}{2\delta y}(u_{i,j+1} -u_{i,j-1}) + \mathcal{O}(\delta y^2)~, \\
    u_{xx}(x_i, y_j, t) &= \frac{1}{\delta x^2}(u_{i+1,j} - 2u_{i,j} +u_{i-1,j}) + \mathcal{O}(\delta x^2)~,\\
    u_{yy}(x_i, y_j, t) &= \frac{1}{\delta y^2}(u_{i,j+1} - 2u_{i,j} +u_{i,j-1}) + \mathcal{O}(\delta y^2)~.
\end{align}
To this end, the explicit time-stepping scheme to update next step solution $u^{n+1}$ is given by:
\begin{align}
   u_{i,j}^n\approx  U_{i,j}^{n+1} &= U_{i,j}^{n}+\delta t f(x_{i},y_{j}, U_{i,j}^{n},U_{i,j-1}^{n},U_{i,j+1}^{n},U_{i+1,j}^{n},U_{i-1,j}^{n} )~, 
\label{eq225}\\
    & \equiv \mathbf{F}(x_{i},y_{j},\delta x, \delta y, \delta t, U_{i,j}^{n},U_{i,j-1}^{n},U_{i,j+1}^{n},U_{i+1,j}^{n},U_{i-1,j}^{n})~,
    \label{eq226}
\end{align}
where $U_{i,j}^{n}$ is the numerical solution at grid point $(x_{i},y_{j}, t_{n})$.

\begin{wrapfigure}{hr}{0.5\textwidth}
  \begin{center}
    \includegraphics[width=0.4\textwidth]{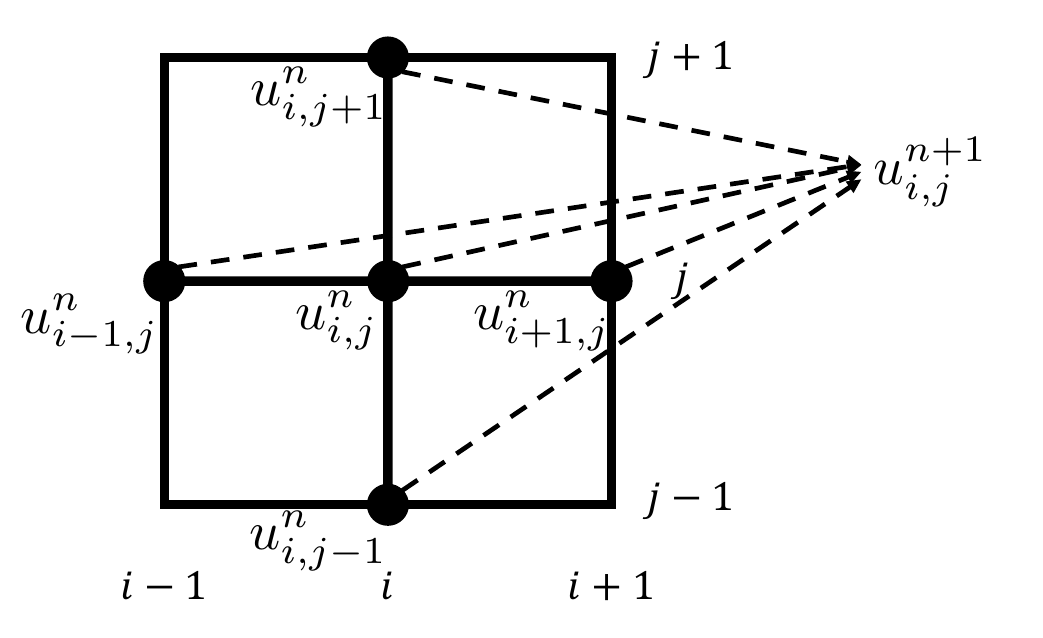}
  \end{center}
  \caption{updating scheme for central difference method}
  \label{fig:fdm}
\end{wrapfigure}

Apparently, the finite difference method (\ref{eq225})  for updating $u^{n+1}$ on a grid point relies on the  previous time steps' solutions, supported on the grid point and its neighbours. The scheme (\ref{eq225})  updates $u^{n+1}_{i,j}$ using five points of $u^{n}$ values (see Figure \ref{fig:fdm}).

Similarly, the finite element method (FEM) approximates the new solution by calculating the corresponded mesh cell coefficient, which is updated by its related nearby coefficients on the mesh.

From this perspective, one may regard the numerical schemes for solving time-dependent PDEs as methods catching the information from  neighbourhood data of interest.

\subsection{Finite Element Method}
\label{sec:fem}

Finite element method (FEM) is a powerful numerical method in solving PDEs.
Consider a 1D wave equation of $u(x,t)$:
\begin{align}
    &u_{tt} - v^2 u_{xx} = f,  \quad x \in [a,b] \equiv \Omega \subset \mathbb{R}, \quad t \in \mathbb{R}^{+} \cup\{0\}~, \\
    &u_x(a, t) =  u_x(b,t) =  0~.
\end{align}
The function $u$ is approximated by a FEM function $u_{h}$ :
\begin{align}
    u(x,t) \approx u_{h}(x,t) &= \sum_{i=1}^{N} a_{i}(t) \psi_{i}(x) \\
\end{align}
where $\psi_{i} \in V$ is the basis functions of some FEM space $V$, and $a_{i}^{n}$ denotes the coefficients.  $N$ denotes the degrees of freedom.\\
Multiply the equation with an arbitrary test function $\psi_{j}$ and integral over the whole domain we have:
\begin{align}
    \int_{\Omega} u_{tt} \psi_{j} \ dx + v^2 \int_{\Omega}\nabla u \nabla \psi_{j} \ dx = \int_{\Omega} f \psi_{j} \ dx \\
\end{align}
and approximate $u(x,t)$ by $u_{h}$:
\begin{align}
     \sum_{i}^{N}\frac{\partial^2 a_{i}(t)}{\partial t^2}\underbrace{\int_{\Omega} \psi_{i} \psi_{j} \ dx}_{\mathbf{M}_{i,j}}  + v^2\sum_{i}^{N}a_{i}(t) \underbrace{\int_{\Omega}\nabla \psi_{i} \nabla \psi_{j} \ dx}_{\mathbf{A}_{i,j}} = \underbrace{\int_{\Omega} f \psi_{j} }_{\mathbf{b}}\ dx ~, \\ 
     \equiv \mathbf{M}^{T}\mathbf{a}_{tt} + v^2 \mathbf{A}^{T} \mathbf{a} = \mathbf{b}~.
\end{align}
Here $\mathbf{M}$ is the mass matrix and $\mathbf{A}$ is the stiffness matrix, $\mathbf{a} =(a_1,..,a_N)^t $ is a $ N$ vector of the coefficients at time $t$. The central difference method for time discretization indicates that \citep{johnson2012numerical}:
\begin{align}
    \mathbf{a}^{n+1} &=2\mathbf{a}^{n} - \mathbf{a}^{n-1} + \mathbf{M}^{-1}(\mathbf{b}-v^2 \mathbf{A}^{T} \mathbf{a}^{n}) ~. 
\label{eq:fem_demo}
\end{align}
This leads to
\begin{align}
    u^{n+1} &\approx u_{h}^{n+1} = \sum_{i}^{N} {a}_{i}^{n+1}\psi_{i}(x)~. \label{eq:fem_demo_2}
\end{align}

\subsection{Long Short-Term Memory}
Long Short-Term Memory networks (LSTM) \citep{hochreiter1997long,graves2005framewise} are a class of artificial recurrent neural network (RNN) architecture that is commonly used for processing sequence data, and can overcome the gradient vanishing issue in RNN. Similar to most RNNs \citep{mikolov2011extensions}, LSTM takes a sequence $\{\vx_1, \vx_2, \cdots, \vx_t\}$ as input and learns hidden vectors $\{\vh_1, \vh_2, \cdots, \vh_t\}$ for each corresponding input. In order to better retain long distance information, LSTM cells are specifically designed to update the hidden vectors. 
The computation process of the forward pass for each LSTM cell is defined as follows:
\begin{align*}
    \vi_t &= \sigma(\mathbf{W}^{(x)}_i \vx_t + \mathbf{W}^{(h)}_i \vh_{t-1} + \mathbf{W}^{(c)}_i \vc_{t-1} + \vb_i)~, \\
    \vf_t &= \sigma(\mathbf{W}^{(x)}_f \vx_t + \mathbf{W}^{(h)}_f \vh_{t-1} + \mathbf{W}^{(c)}_f \vc_{t-1} + \vb_f)~, \\
    \vc_t &= \vf_t \vc_{t-1} + \vi_t \tanh (\mathbf{W}^{(x)}_c \vx_t + \mathbf{W}^{(h)}_c \vh_{t-1} + \vb_c)~, \\
    \vo_t &= \sigma(\mathbf{W}^{(x)}_o \vx_t + \mathbf{W}^{(h)}_o \vh_{t-1} +\mathbf{W}^{(c)}_o \vc_t + \vb_o), \\
    \vh_t &= \vo_t \tanh (\vc_t)~,
\end{align*}
where $\sigma$ is the logistic sigmoid function, $\mathbf{W}$s are weight matrices, $\vb$s are bias vectors, and subscripts $\vi$, $\vf$, $\vo$ and $\vc$ denote the input gate, forget gate, output gate and cell vectors respectively, all of which have the same size as hidden vector $\vh$.

This  LSTM structure is used in the paper to simulate the numerical solutions of partial differential equations.

\section{Proposed Method}
\label{sec:method}
\subsection{Mathematical Motivation}
\label{sec:math}
Recurrent neural network including LSTM is an artificial neural network structure of the form~\citep{lipton2015critical}:
\begin{equation}
    \vh^t =\sigma (\rmW^{hx}\vx^t+\rmW^{hh}\vh^{t-1}+\vb_h)
    \equiv
    \sigma_a(\vx^t,\vh^{t-1}) \equiv \sigma_b(\vx^{0},\vx^{1},\vx^2,\cdots,\vx^t)~,
    \label{eq311}
\end{equation}
where $\vx^t \in \mathbb{R}^d$ is the input data of the $t^{th}$ state and $\vh^{t-1} \in \mathbb{R}^h$ denotes the processed value in its previous state by the hidden layers. The  output $\vy^{t}$ of the current state  is updated by the current state value $\vh^{t}$:
\begin{align}
    \vy^t &= \sigma (\rmW^{hy} \vh^{t} +
    \vb_y)\\
    &\equiv
    \sigma_c(\vh^t)\equiv \sigma_d(\vx^{0},\vx^{1},\vx^2,\cdots,\vx^t)~.
\end{align}

Here $\rmW^{hx} \in \mathbb{R}^{h\times d}$, $\rmW^{hh}\in \mathbb{R}^{h\times h}$, $\rmW^{hy}\in \mathbb{R}^{h\times h}$ are the matrix of weights,  vectors $\vb_h,\vb_y \in \mathbb{R}^h$ are the coefficients of bias, and $\sigma, \sigma_a, \sigma_b, \sigma_c, \sigma_d$ are corresponded activation and mapping functions.
With proper design of input and forget gate, LSTM can effectively yield a better control over the gradient flow and better preserve useful information from long-range dependencies \citep{graves2005framewise}.

Now consider a temporally continuous vector function $\vu \in \mathbb{R}^n$ given by an ordinary differential equation with the form:
\begin{equation}
    \frac{d \vu(t)}{d t} = g(\vu(t))~.
\end{equation}

Let $\vu^{n} =\vu(t = n\delta t)$, a forward Euler's method for solving $\vu$ can be easily derived from the Taylor's theorem which  gives the following first-order accurate approximation of the time derivative:
\begin{equation}
    \frac{d \vu^n}{d t} = \frac{\vu^{n+1} -\vu^{n}}{\delta t} + \mathcal{O}(\delta t)~. \label{eq313}
\end{equation}
Then we  have:
\begin{align}
    \frac{d \vu}{d t} =  g(\vu) \xrightarrow{(\ref{eq313})} \vu^{n+1} &=  \vu^{n} + \delta t \ g(\vu^{n}) + \mathcal{O}(\delta t^2)\nonumber \\
    &\rightarrow \hat{\vu}^{n+1} =
    f_1(\hat{\vu}^{n}) =\underbrace{f_1\circ f_1\circ\cdots f_1(\hat{\vu}^{0})}_{n} 
    \label{eq315}
\end{align}
Here  $\hat{\vu}^n \approx \vu(n\delta t)$ is the numerical approximation and $f_1 \equiv  \vu^{n} + \delta t \ g(\vu^{n}) :\mathbb{R}^n \rightarrow \mathbb{R}^n$. Combining equations (\ref{eq311}) and (\ref{eq315}) one may notice that the residual networks, recurrent neural network and also LSTM networks can be regarded as a numerical scheme for solving time-dependent differential equations if more layers are added and  smaller time steps are taken.~\citep{chen2018neural}

Canonical structure for such recurrent neural network usually calculates the current state value by its previous time step value $\vh^{t-1}$ and current state input $\vx^t$. Similarly,  in numerical PDEs, the next step data at a grid point is updated from the previous (and current) values on its nearby grid points (see Eq.~\ref{eq225}). 

Thus, what if we replace the temporal input $\vh^{t-1}$ and $\vx^{t}$ with spatial information? A  simple sketch of the upwinding method for a $1D$ example of $u(x,t)$:
\begin{equation}
   u_t + \nu u_{x} = 0 \label{eq316}
\end{equation}
will be:
\begin{align}
    u_{i}^{n+1} & = u_{i}^{n} -\nu\frac{\delta t}{\delta x} (u_{i}^n - u_{i-1}^n)   + \mathcal{O}(\delta x,\delta t) 
    \rightarrow \hat{u}_{i}^{n+1} = f_2(\hat{u}_{i-1}^{n},\hat{u}_{i}^{n})  \\
      &\equiv f_\theta \big(f_\eta(\vx_{i},\vh_{i-1}(u))\big)
    = f_{\theta,\eta}\big(\hat{u}_{0}^n, \hat{u}_{1}^n, \cdots, \hat{u}_{i-1}^n,\hat{u}_{i}^n\big)=v_{i}^{n+1} \\
    &\vx_{i} = \hat{u}_{i}^{n}, \;
    \vh_{i-1}(\hat{u}) = \sigma ( \hat{u}_{i-1}^n,\vh_{i-2}(\hat{u})) \equiv f_\eta(\hat{u}_{0}^{n}, \hat{u}_{1}^{n}, \hat{u}_{2}^{n},\cdots, \hat{u}_{i-1}^{n}).
\end{align}
Here we use $v_{i}^{n+1}$ to denote the prediction of $\hat{u}_{i}^{n+1}$ processed by neural network. We replace the temporal previous state $\vh^{t-1}$with spacial grid value $\vh_{i-1}$ and input the numerical solution $\hat{u}_i^{n}\approx u(i\delta x,  n\delta t)$ as current state value, which indicates the neural network could be seen as a forward Euler method for equation \ref{eq316}~\citep{lu2018beyond}. Function $f_2 \equiv  \hat{u}_i^{n} - \nu  \frac{\delta t}{\delta x}(\hat{u}_i^n-\hat{u}_{i-1}^n)  :\mathbb{R}^2 \rightarrow \mathbb{R}$ and the function $f_\theta$ represents the dynamics of the hidden layers in decoder with parameters $\theta$, and $f_\eta$ specifies the dynamics of the LSTM layer~\citep{hochreiter1997long,graves2005framewise}  in encoder withe parameters $\eta$. The function $f_{\theta, \eta}$ simulates the dynamics of the  Neural-PDE with paramaters $\theta$ and $\eta$.
By applying Bidirectional neural network, all grid data are transferred and it enables LSTM to simulate the  PDEs as :
\begin{align}
    v_{i}^{n+1} &= f_\theta\big(f_\eta(\vh_{i+1}(\hat{\hat{u}}), \hat{u}_{i}^{n},\vh_{i-1}(\hat{u}))\big) \\
    &\vh_{i+1}(\hat{u}) \equiv f_\eta(\hat{u}_{i+1}^{n}, \hat{u}_{i+2}^{n}, \hat{u}_{i+3}^{n},\cdots, \hat{u}_{k}^{n}).
\end{align}
For a time-dependent PDE, if we map all our  grid data into an input matrix which contains the information of $\delta x, \delta t$, then the neural network would regress such coefficients as constants and will learn and filter the physical rules from all the $k$ mesh grids data as:
\begin{equation}
    v_{i}^{n+1} = f_{\theta, \eta}\big( \hat{u}_{0}^{n}, \hat{u}_{1}^{n}, \hat{u}_{2}^{n},\cdots, \hat{u}_{k}^{n} \big)
\end{equation}

The LSTM neural network is designed to overcome the vanishing gradient issue through hidden layers, therefore we use such recurrent structure to increase the stability of the numerical approach in deep learning. The highly nonlinear function $f_{\theta,\eta}$ simulates the dynamics of  updating rules for $u_{i}^{n+1}$, which works in a way  similar to a finite difference method (section \ref{sec:fdm}) or a finite element method.

\subsection{Neural-PDE}
\begin{wrapfigure}{hr}{0.38\textwidth}
  \begin{center}
    \includegraphics[width=0.35\textwidth]{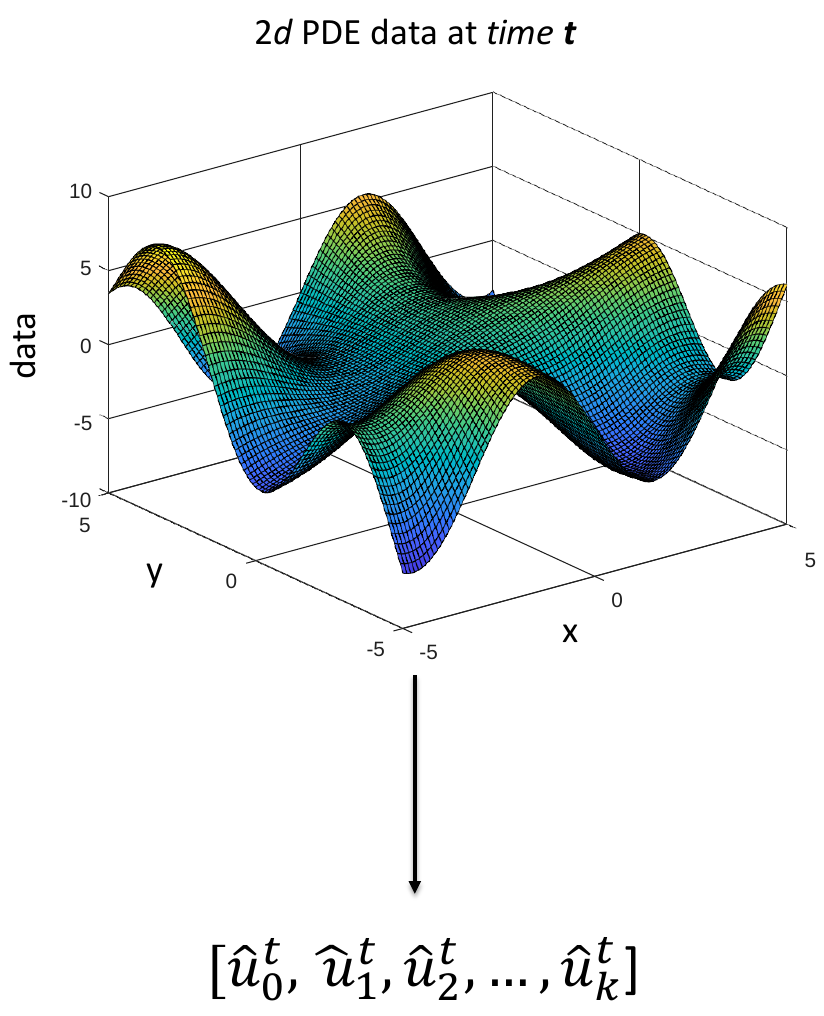}
  \end{center}
  \caption{An example of maping $2d$ data matrix into $1d$ vector  where $k =N_{x} \times N_{y}$  and $N_{x}$ and $N_{y}$ are the numbers of grid points on $x$ and $y$, respectively.}
  \label{fig:2d_1d}
\end{wrapfigure}

In particular, we use the bidirectional LSTM \citep{hochreiter1997long,graves2005framewise} to better retain the state information from data on grid points which are neighbourhoods in the mesh but far away in input matrix. 

\begin{figure}[t]
\centering
\includegraphics[width=1.1\linewidth]{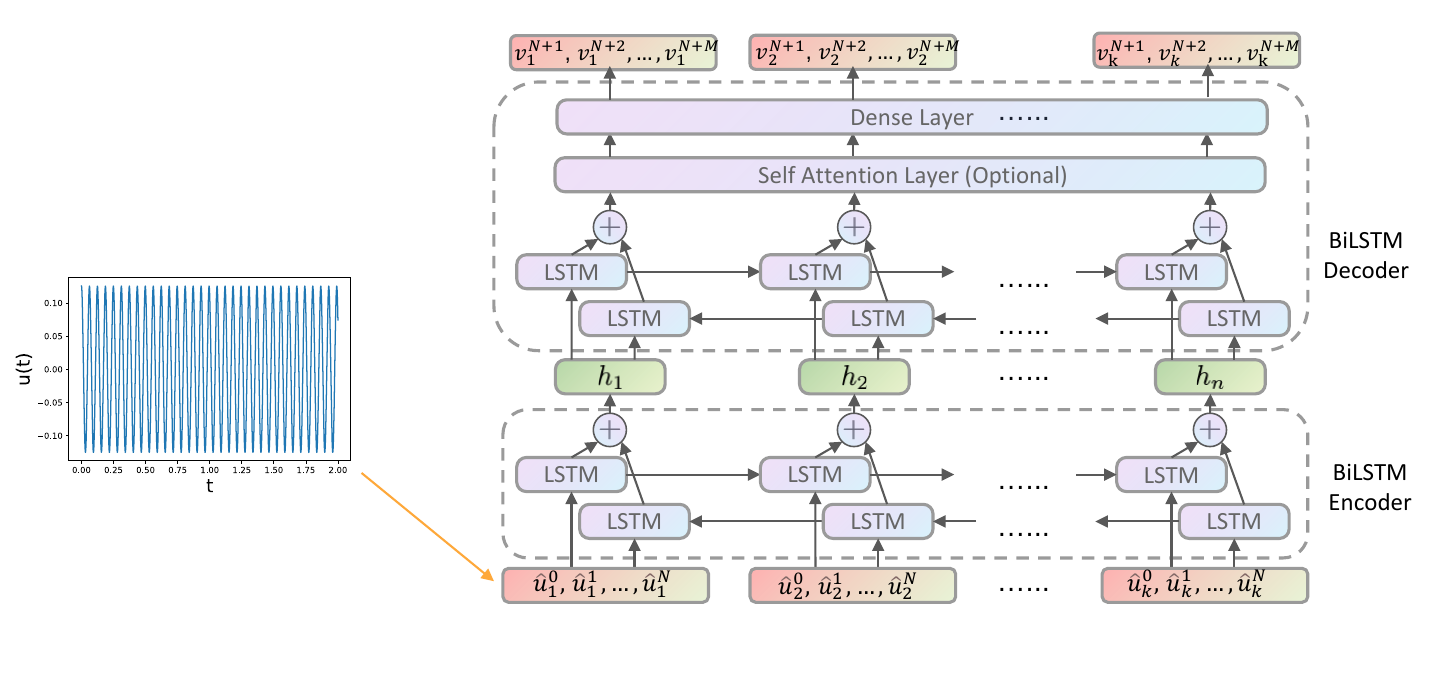}  
\caption{Neu-PDE}
\label{fig:model}
\end{figure}

The right frame of Figure~\ref{fig:model} shows the overall design of the Neural-PDE. Denote the time series data at collocation points as $\va_1^N, \va_2^N, \cdots, \va_k^N$  with $\va_{i}^N=[\hat{u}_i^0, \hat{u}_i^1, \cdots, \hat{u}_i^N]$ at $i^{th}$ point. The superscript represents different time points. The Neural-PDE takes the past states $\{\va_1^N, \va_2^N, \cdots, \va_k^N\}$ of all collocation points,  and outputs the predicted future states $\{\vb_1^M, \vb_2^M, \cdots, \vb_k^M\}$, where $\vb_i^M= [v_i^{N+1}, v_i^{N+2}, \cdots, v_i^{N+M}]$ is the Neural-PDE prediction for the $i^{th}$ collocation point at time points from $N+1$ to $N+M$. The data from time point 0 to $N$ are the  training data set. 

The Neural-PDE is an encoder-decoder style sequence model that first maps the input data to a low dimensional latent space that
\begin{equation}
    \vh_i = \overrightarrow{\mathrm{LSTM}}(\va_{i}) \oplus \overleftarrow{\mathrm{LSTM}}(\va_{i}),
\end{equation}
where $\oplus$ denotes concatenation and $\vh_i$ is the latent embedding of point $\va_i$ under the environment.

One then decoder, another bi-lstm with a dense layer:
\begin{equation}
    v_i = \left( \overrightarrow{\mathrm{LSTM}}(\vh_i) \oplus \overleftarrow{\mathrm{LSTM}}(\vh_i) \right) \cdot \mathbf{W},
\end{equation}
where $\mathbf{W}$ is the learnable weight matrix in the dense layer. Moreover, the final decode layers could also have an optional self attention~\citep{vaswani2017attention} layer, which makes the model easier to learn long-range dependencies of the mesh grids.

During training process, mean squared error (MSE) loss $\mathcal{L}$ is used as we typically don't know the specific form of the PDE.
\begin{equation}
    \mathcal{L} = \sum_{t=N+1}^{N+M}\sum_{i=1}^{k} ||\hat{u}_{i}^{t} - v_{i}^{t} ||^2~,
\end{equation}

\subsection{Data Initialization and Grid Point Reshape}

In order to feed the  data into our sequence model framework, we  map the PDE solution data onto a $K \times N$ matrix, where $K \in \mathbb{Z}^+$ is the dimension of the grid points and  $N \in \mathbb{Z}^+$ is the length of the time series data on each grid point. There is no regularization for the input order of the grid points data in the matrix because of the bi-directional structure of the Neural-PDE.  For example, a $2D$ heat equation at some time $t$ is reshaped into a $1D$ vector (See Fig.~ \ref{fig:2d_1d}). Then the matrix is formed accordingly.

For a $n$-dimensional time-dependent partial differential equation with $K$ collocation points, the input and output data  for $t \in (0,T)$ will be of the form:
\begin{align}\label{eq:3.3.1}
&\scale[0.85]{\mA(K,N) =
\begin{bmatrix}
\va_{0}^N \\
\vdots \\
\va_{\ell}^N \\
\vdots \\
\va_{K}^N \\
\end{bmatrix}
=
\begin{bmatrix}
\hat{u}_{0}^0 & \hat{u}_{0}^1 & \cdots & \hat{u}_{0}^n & \cdots &\hat{u}_{0}^N\\
\vdots & \vdots & \ddots & \vdots& \ddots & \vdots \\
\hat{u}_{\ell}^0 & \hat{u}_{\ell}^1 & \cdots & \hat{u}_{\ell}^n & \cdots &\hat{u}_{\ell}^N\\
\vdots & \vdots & \ddots & \vdots & \ddots & \vdots \\
\hat{u}_{K}^0 & \hat{u}_{K}^1 & \cdots & \hat{u}_{K}^n & \cdots &\hat{u}_{K}^N
\end{bmatrix}}\\
&\scale[0.85]{\mB(K,M)  = 
    \begin{bmatrix}
    \vb_{0}^M \\
    \vdots \\
    \vb_{\ell}^M \\
    \vdots \\
    \vb_{K}^M \\
    \end{bmatrix}
    =
\begin{bmatrix}
v_{0}^{N+1} & v_{0}^{N+2} & \cdots & v_{0}^{N+m} & \cdots &v_{0}^{N+M}\\
\vdots & \vdots & \ddots & \vdots & \ddots & \vdots \\
v_{\ell}^{N+1} & v_{\ell}^{N+2} & \cdots & v_{\ell}^{N+m} & \cdots &v_{k}^{N+M}\\
\vdots & \vdots & \ddots &\vdots & \ddots & \vdots \\
v_{K}^{N+1} & v_{K}^{N+2} & \cdots & v_{K}^{N+m} & \cdots &v_{K}^{N+M}\\
\end{bmatrix}}
\end{align}

Here $N = \frac{T}{\delta t}$ and each row $\ell$ represents the time series data at the  $\ell^{th}$ mesh grid, and  $M$ is the time length of the predicted data.

By adding Bidirectional LSTM  encoder in the Neural-PDE, it will automatically extract the information from the time series data  as:
\begin{equation}
    \mB(K,M)  
    = PDESolver(\mA(K,N)) =PDESolver(\va_{0}^N,\va_{1}^N, \cdots \; \va_{i}^N, \cdots \;, \va_{K}^N)
\end{equation}

\section{Computer Experiments}
\label{sec:experiments}

\begin{table}[htbp]
\centering
\sisetup{round-mode=places
        ,round-precision=3}
\caption{Error analysis models} 
\begin{tabular}{lccc}
\toprule
& Wave & Heat & Burgers'\\
\midrule
Equation      & $u_{tt} = \frac{1}{16\pi^2}u_{xx}$    & $u_{t} = u_{xx}$ & $\frac{\partial u}{\partial t}+u \frac{\partial u}{\partial x}=0.1 \frac{\partial^{2} u}{\partial x^{2}}$   \\
IC      & $\sin(4\pi x)$    & $6\sin(\pi x)$ & $u(0 \leq x \leq L, t=0) = 0.9$   \\
BC      & $periodic$    & $periodic$ & $periodic$   \\
\bottomrule
\end{tabular}
\label{tab:1}
\end{table}

\begin{table}[htbp]
\centering
\sisetup{round-mode=places
        ,round-precision=3}
\caption{$L^2$ error for model evaluation.} 
\begin{tabular}{lccc}
\toprule
$\Delta x = 0.1$ & Wave & Heat & Burgers'\\
\midrule
$\Delta t = 0.1$      & \num{4.385e-03}    & \num{6.912e-5} & \num{9.450e-04}   \\
$\Delta t = 0.01$      & \num{3.351e-05}    & \num{5.809e-5} & \num{5.374e-3}   \\
$\Delta t = 0.001$      & \num{1.311e-05}    & \num{3.757e-5} & \num{1.244e-3}   \\
\bottomrule
\end{tabular}
\label{tab:2}
\end{table}

\begin{table}[htbp]
\centering
\sisetup{round-mode=places
        ,round-precision=3}
\caption{$L^2$ error for model evaluation.} 
\begin{tabular}{lccc}
\toprule
$\Delta t = 0.1$ & Wave & Heat & Burgers'\\
\midrule
$\Delta x = 0.1$      & \num{2.190e-05}    & \num{1.162e-4} & \num{2.561e-4}   \\
$\Delta x = 0.01$      & \num{6.059e-05}    & \num{7.706e-4} & \num{4.206e-4}   \\
$\Delta x = 0.001$      & \num{1.498e-05}    & \num{1.400e-5} & \num{3.700e-4}   \\
\bottomrule
\end{tabular}
\label{tab:3}
\end{table}

Since the Neural-PDE is a sequence to sequence learning framework which allows to predict within any time period by the given data. One may test the Neural-PDE using different permutations of training and predicting time periods for its efficiency, robustness and accuracy. In the following examples, the whole dataset is randomly splitted  in $80\%$ for training and $20\%$ for testing. We will predict the next $t_p \in [31 \times\delta t, 40 \times \delta t$] PDE solution by using its previous $t_{tr}\in [0, 30 \times\delta t]$ data as:
\begin{equation}
    \mB(K,10) = PDESolver(\mA(K,30))
\end{equation}

We tested Neu-PDE using three classical PDE models with different $\Delta x$ and $\Delta t$, Table ~\ref{tab:1} summarizes the information of these models.
Table~\ref{tab:2} and Table~\ref{tab:3} show the experimental results of the Neural-PDE model solving the above three different PDEs.
We used the Neural-PDE which only consists of 3 layers: 2 bi-lstm (encoder-decoder) layers with 20 neurons each and 1 dense output layer with 10 neurons and achieved MSEs from $ \mathcal{O}( 10^{-3})$ to $\mathcal{O} (10^{-5})$ within 20 epochs, a MLP based neural network such as Physical Informed Neural Network~\citep{raissi2019physics} usually will have more layers and neurons to achieve similar $L^2$ errors. Additional examples are also discussed in this section.

\subsection*{Example: Wave equation}
Consider the $1D$ wave equation:
\begin{align}
    & u_{tt} = cu_{xx}, \; x \in [  0, 1  ], \; t \in [0,2]~, \\
    & u(x,0) = sin(4\pi x) \\
    & u(0, t) = u(1, t) 
\end{align}
Let $c = \frac{1}{16\pi^2}$ and use the analytical solution given by the characteristics for the training and testing data:

\begin{equation}
    u(x,t) = \frac{1}{2}(sin(4\pi x +t) +sin(4\pi x -t))~.
\end{equation}

\begin{figure}[H]
\centering  
\subfigure[Exact  Test Dataset]{\includegraphics[width=0.45\linewidth]{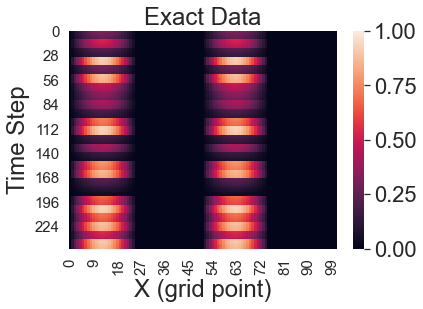}}
\subfigure[Predicted  Test Dataset]{\includegraphics[width=0.45\linewidth]{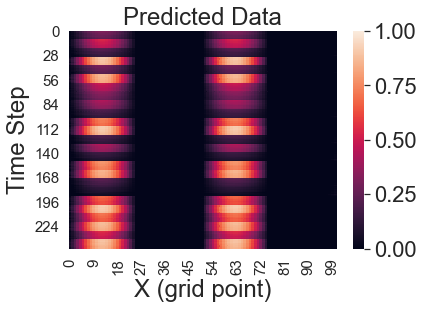}}
\subfigure[Training Metrics]{\includegraphics[width=.4\linewidth]{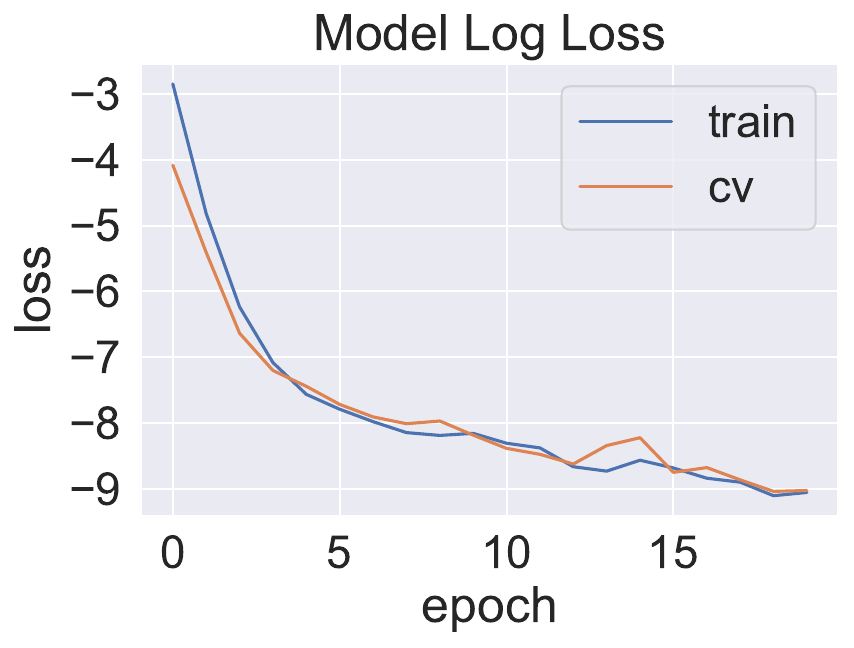}}
\caption{The Neural-PDE for solving the wave equation.}
\medskip
\small
\label{fig:wave}
\end{figure}

Here we used $\delta x = 1\times 10^{-2}$, $\delta t = 1\times 10^{-2}$, and the mesh grid size is $101$. We obtained a MSE $3.5401\times 10^{-5}$.  The test dataset batch size is $25$ and thus the total discrete testing time period is $250$. Figures~\ref{fig:wave}$(a)$ and \ref{fig:wave}$(b)$ are the heat map for the exact test data and our predicted test data. Figure~\ref{fig:wave}$(c)$ shows both training and cross-validation errors of Neural-PDE convergent within $20$ epochs.

\begin{figure}[H]
\centering  
\subfigure[$t=1.991$]{\includegraphics[width=0.4\linewidth]{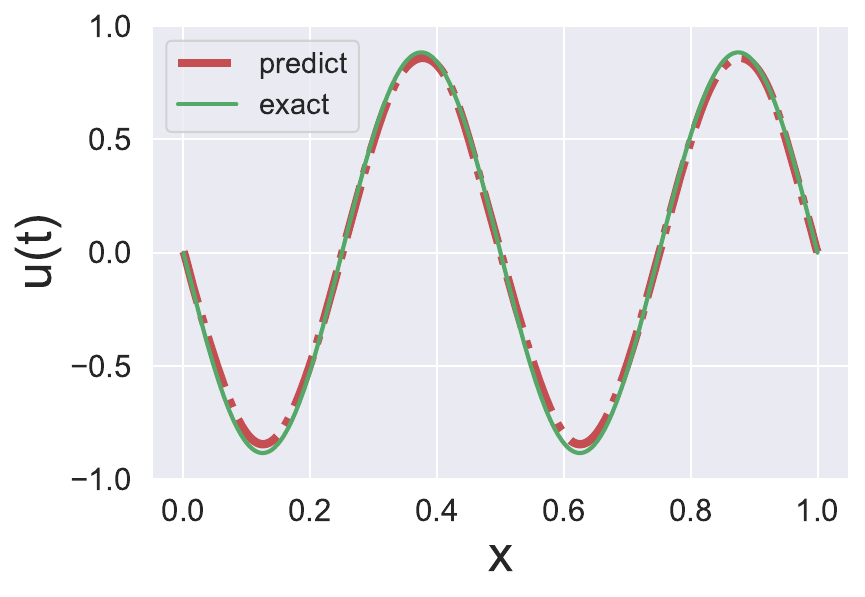}}
\subfigure[$t=1.994$]{\includegraphics[width=0.4\linewidth]{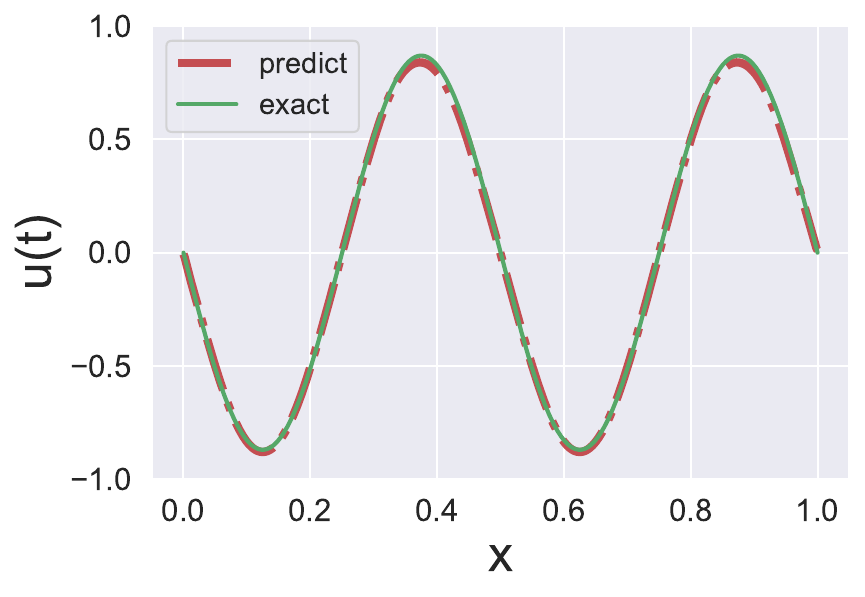}}
\subfigure[$t=1.997$]{\includegraphics[width=0.4\linewidth]{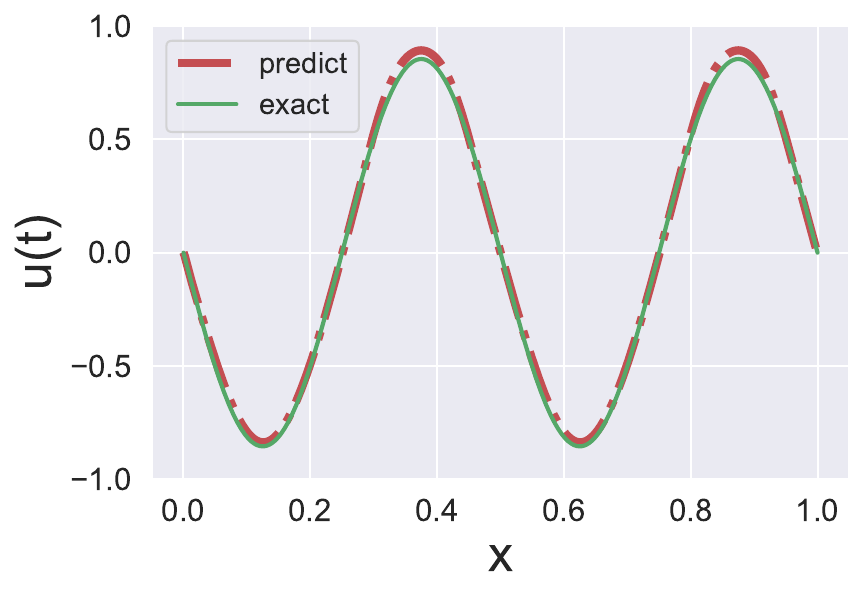}}
\subfigure[$t=2$]{\includegraphics[width=0.4\linewidth]{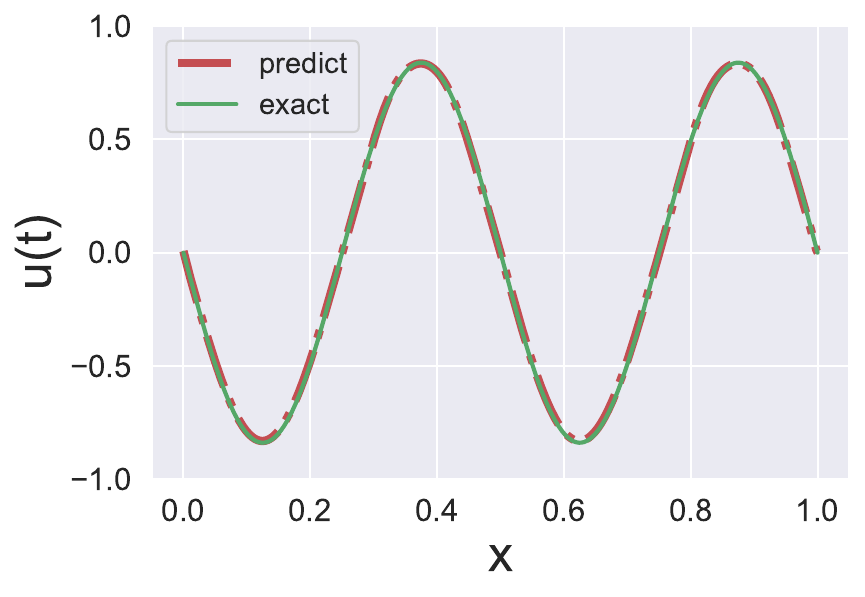}}
\caption{Comparison between exact solution and Neural-PDE prediction of the 1D wave equation at various time points.}
\medskip
\small
\label{fig:wave_final}
\end{figure}

We selected the final four states for computation and compared them with analytic solutions. The result indicates that the Neural-PDE is robust in capturing the physical laws of wave equation and predicting the sequence time period. See Figure~\ref{fig:wave_final}.

\subsection*{Example: Heat equation}
The 1D wave equation case maps the data into a matrix (\ref{eq:3.3.1}) with its original spatial locations. In this test, we solve  the 2D heat equation describing how the motion or diffusion of a heat flow evolves over time. Here the  2-dimensional PDE grid in space is mapped into matrix without regularization of the position. The experimental results show that the Neural-PDE is able to capture the valuable features regardless of the order of the  grid points in the matrix. Let's start with a  2D heat equation as follows:
\begin{align}
&u_{t} = u_{xx} +u_{yy}~,\\
  &u(x,y,0)=\left\{
  \begin{array}{@{}ll@{}}
    0.9, & \text{if}\ (x-1)^2 + (y-1)^2<0.25 \\
    0.1, & \text{otherwise}
  \end{array}\right. \\
  &\Omega = [0,2]\times [0,2], \; t \in[0,0.15]~.
\end{align}

\begin{figure}[t]
\centering  
\subfigure[Exact Test Dataset]{\includegraphics[width=0.45\linewidth]{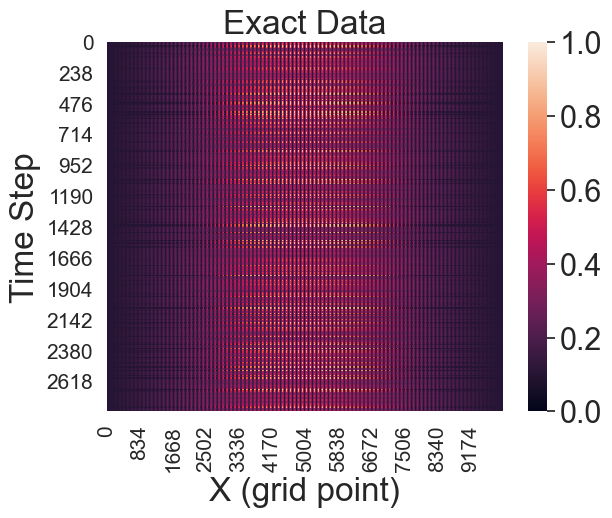}}
\subfigure[Predicted Test Dataset]{\includegraphics[width=0.45\linewidth]{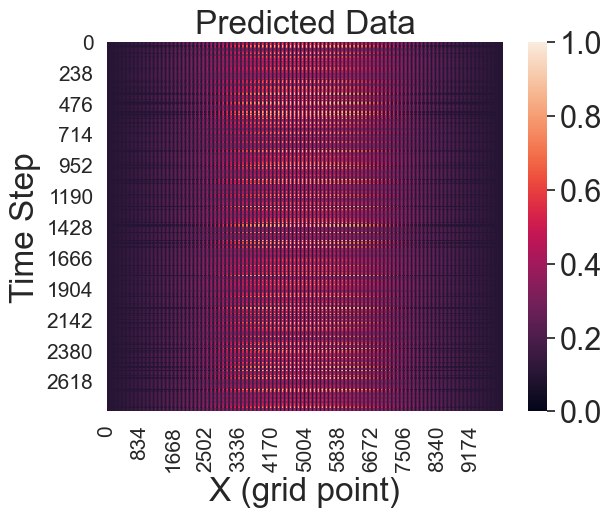}}
\caption{Heatmaps of the heat equation data. $\delta x = 0.02, \delta y = 0.02$, $\delta t = 10^{-4}$, MSE: $2.1551 \times 10^{-6}$.}
\label{fig:heat_eq_heatmap}
\end{figure}

\begin{figure}[H]
\centering  
\subfigure[]{\includegraphics[width=0.4\linewidth]{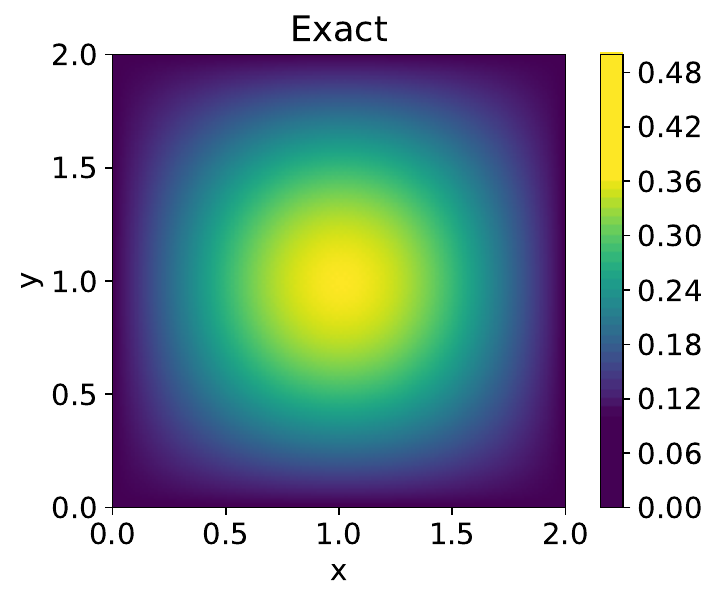}}
\subfigure[]{\includegraphics[width=0.4\linewidth]{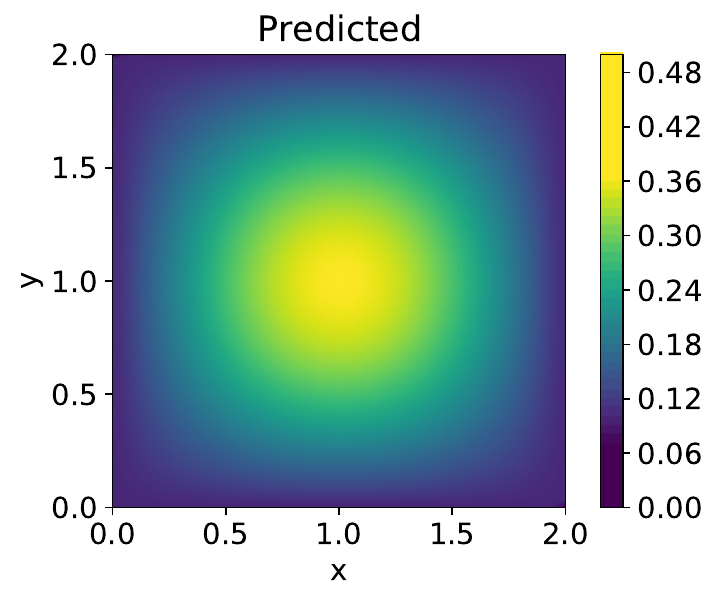}}
\subfigure[]{\includegraphics[width=0.4\linewidth]{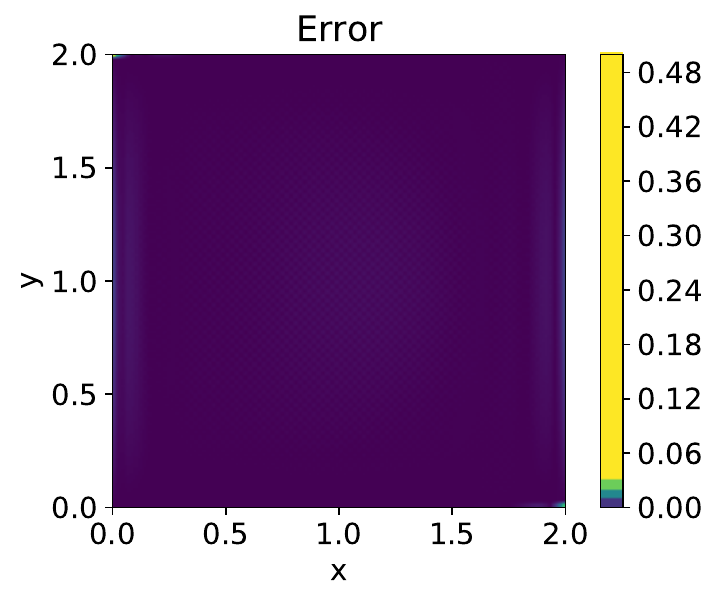}}
\subfigure[Training Metrics]{\includegraphics[width=0.4\linewidth]{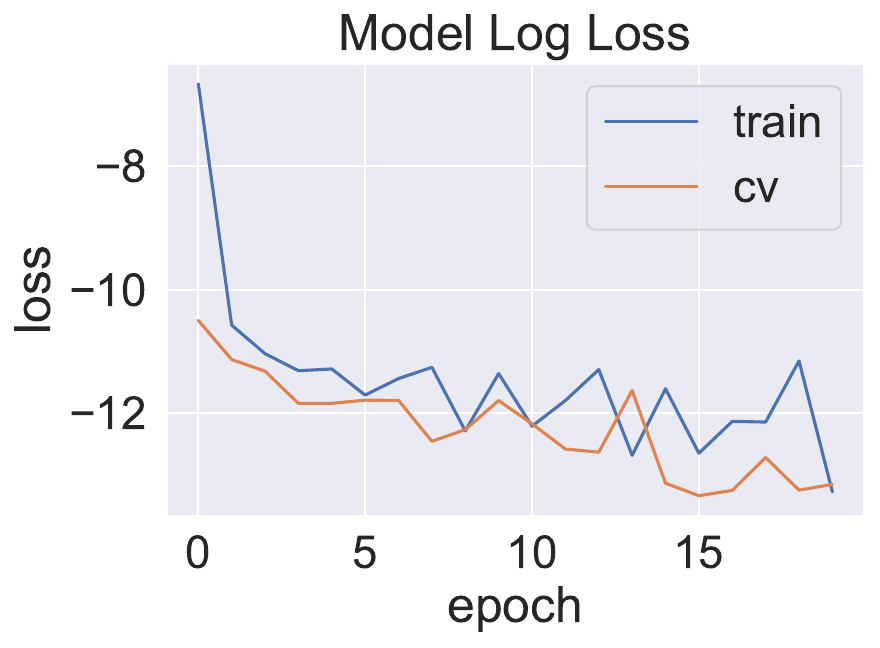}}
\caption{ The Neural-PDE for solving the 2D heat equation. $(a)$ is the exact solution $u(x,y, t = 0.15)$ at the final state. $(b)$ is the Neural-PDE prediction. $(c)$ is the corresponding error map and $(d)$ shows the training and cross-validation errors.}
\label{fig:final_heat}
\end{figure}

Figures~\ref{fig:heat_eq_heatmap} and \ref{fig:final_heat} show the test of the Neural-PDE using the 2D heat equation. We obtained a MSE $\mathcal{O}(10^{-6})$. 

\subsection*{Example: Inviscid Burgers' Equation}
\begin{figure}[H]
\centering  
\subfigure[Exact  Test Dataset]{\includegraphics[height=5cm,width=0.45\linewidth]{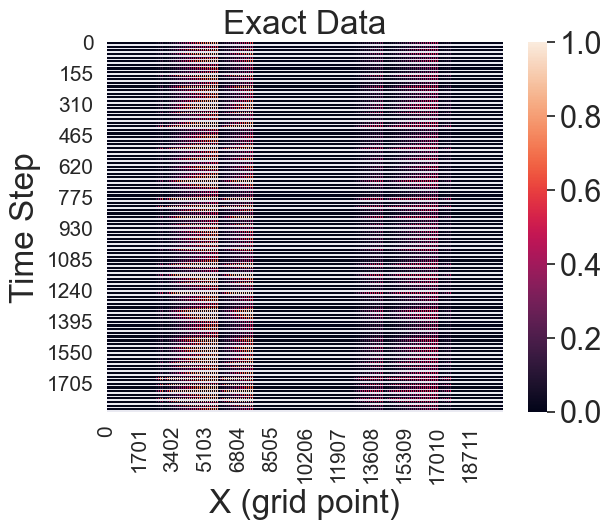}}
\subfigure[Predicted  Test Dataset]{\includegraphics[height=5cm,width=0.45\linewidth]{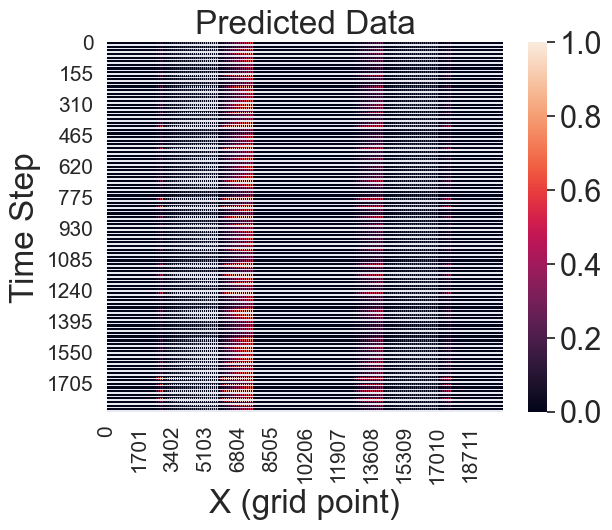}}
\subfigure[Training Metrics]{\includegraphics[width=0.4\linewidth]{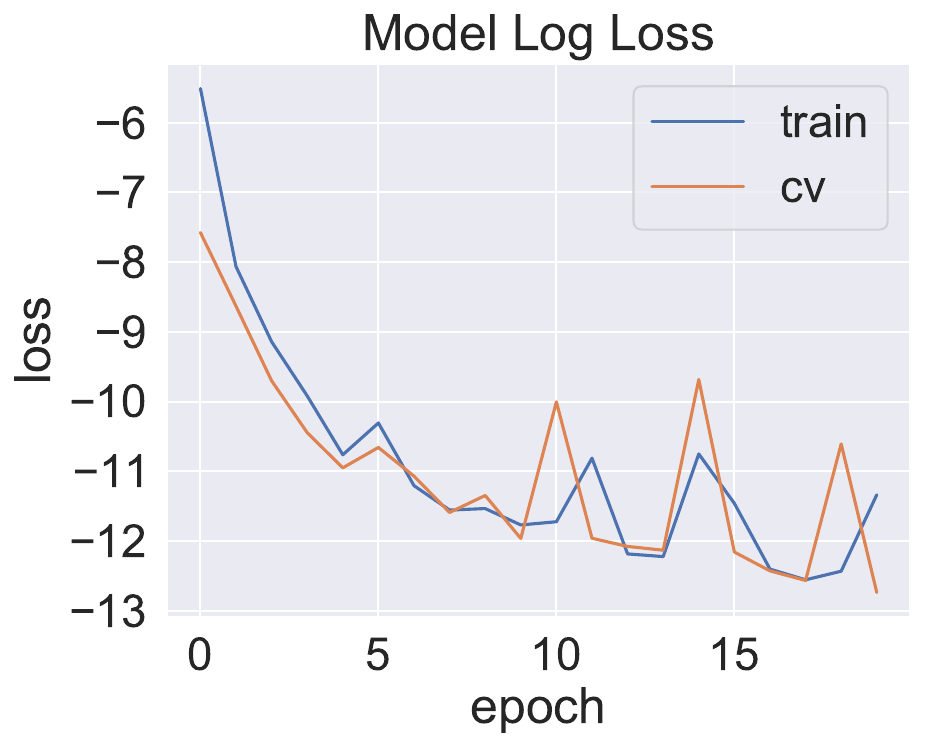}}
\caption{Neural-PDE prediction on the $2D$ Burgers' equation.}
\label{fig:Burgers_heat}
\end{figure}

\begin{figure}[H]
\centering  
\subfigure[Exact $u( t=1)$]{\includegraphics[width=0.45\linewidth]{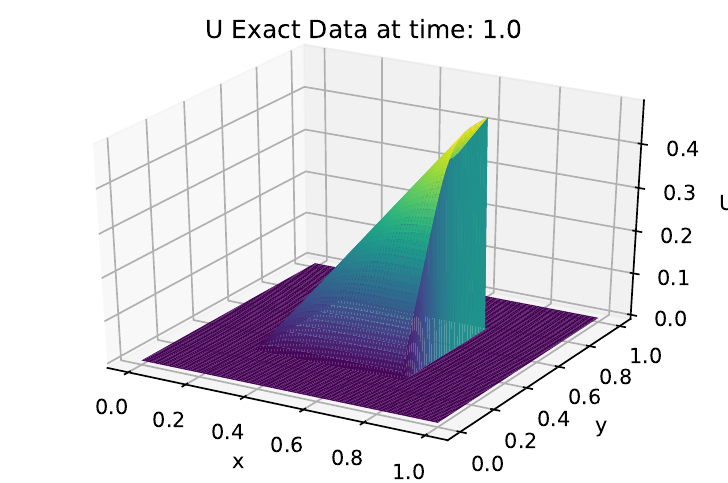}}
\subfigure[Predicted $u(t=1)$]{\includegraphics[width=0.45\linewidth]{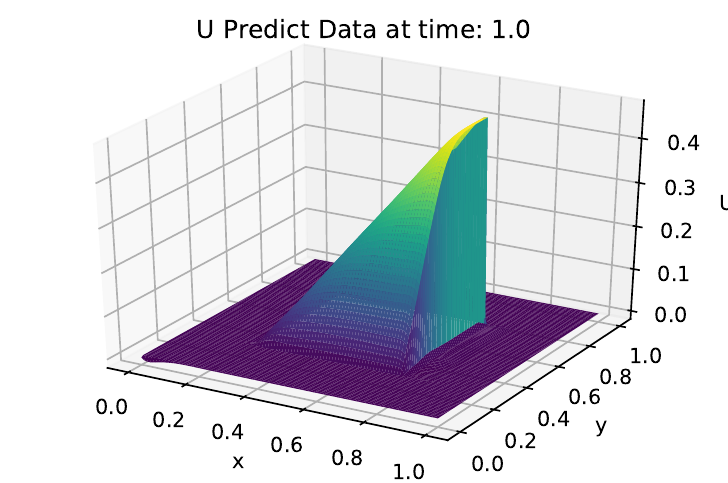}}
\subfigure[Exact $v(t=1)$]{\includegraphics[width=0.45\linewidth]{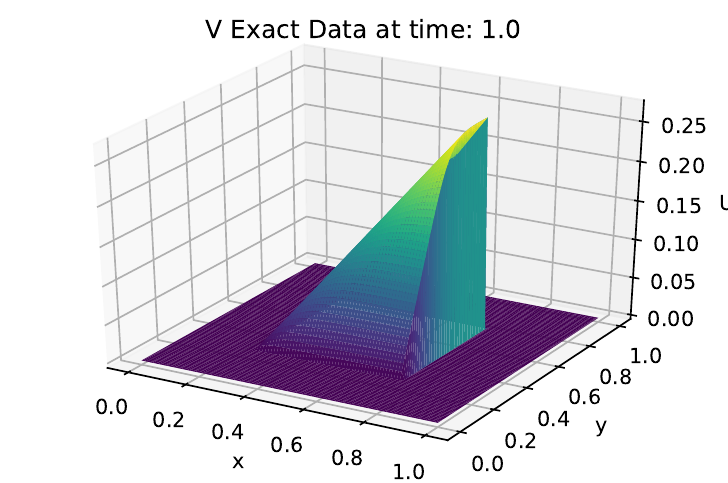}}
\subfigure[Predicted  $v(t=1)$]{\includegraphics[width=0.45\linewidth]{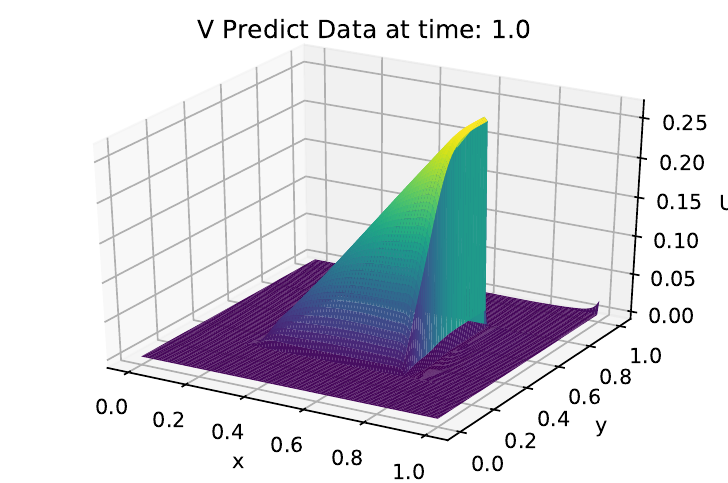}}
\caption{Neural-PDE shows accurate prediction on Burgers' equation.}
\label{fig:Burgers_value}
\end{figure}

Inviscid Burgers' equation is a classical nonlinear PDE in fluid dynamics. In this example, we consider a $2D$ invicid Burgers' equation which has the following hyperbolic  form:
\begin{align}
    &\frac{\partial u}{\partial t} + u\frac{\partial u}{\partial x} +v\frac{\partial u}{\partial y} = 0~, \quad
    \frac{\partial v}{\partial t} + u\frac{\partial u}{\partial x} +v\frac{\partial u}{\partial y} = 0~, \\
    &\Omega = [0,1]\times[0,1], t\in [0, 1]~,
\end{align}
and with the initial and boundary conditions:
\begin{align}
    &u(0.25 \leq x \leq 0.75, \ 0.25 \leq y \leq 0.75, t=0) = 0.9~,\\
    &v(0.25 \leq x \leq 0.75,\ 0.25 \leq y \leq 0.75, t=0) = 0.5~,\\
    &u(0,y,t) = u(1,y,t) = v(x,0,t) =v(x,1,t) =0 ~.
\end{align}
The invicid Burgers' equation is difficult to solve due to the  discontinuities (shock waves) in the solutions. We use a upwinding  finite difference scheme to create the training data and put the velocity $u,v$ in to the input matrix. Let $\delta x = \delta y = 10^{-2}, \delta t = 10^{-3}$, our empirical results (see Figure \ref{fig:Burgers_value}) show that the Neural-PDE is able to learn the shock waves, boundary conditions and the rules of the equation, and predict $u$ and $v$ simultaneously with an overall MSE of $2.3070\times 10^{-6}$. The heat maps of exact solution and predicted solution are shown in Figure \ref{fig:Burgers_heat}.

\subsection*{Example: Multiscale Modeling: Coupled Cahn–Hilliard–Navier–Stokes System}
Finally, let's consider the following $2D$ Cahn–Hilliard–Navier–Stokes system widely used for modeling complex fluids:
\begin{align}
&\vu_{t} +\vu \cdot \nabla \vu
= -\nabla p + \nu \Delta \vu  - \phi  \nabla \mu ~,\label{eq4.1}\\
    &\phi_t{} + \nabla \cdot (\vu\phi) = M\Delta \mu ~,\label{eq4.2} \\
    &\mu = \lambda(-\Delta\phi + \frac{\phi}{\eta^2}( \phi^2-1))~, \label{eq4.3} \\
    & \nabla \cdot \vu = 0 ~.\label{eq4.4}
\end{align}
In this  example we  use the following initial condition:
\begin{align}
    \phi(x,y,0) &= (\frac{1}{2}-50\tanh(f_1-0.1)) + (\frac{1}{2}-50\tanh(f_2-0.1))~,~{\rm where} \\
    f_1(x,y) &= \sqrt{(x+0.12)^2+(y)^2}, \; f_2(x,y) = \sqrt{(x-0.12)^2+(y)^2} ~\\ 
    \text{with}\;  x  &\in [-0.5,0.5], \; y \in [-0.5,0.5], \; t \in [0,1], \; M = 0.1, \; \nu = 0.01, \; \eta = 0.1.
\end{align}

\begin{figure}[t]
\centering  
\subfigure{\includegraphics[width=0.24\linewidth]{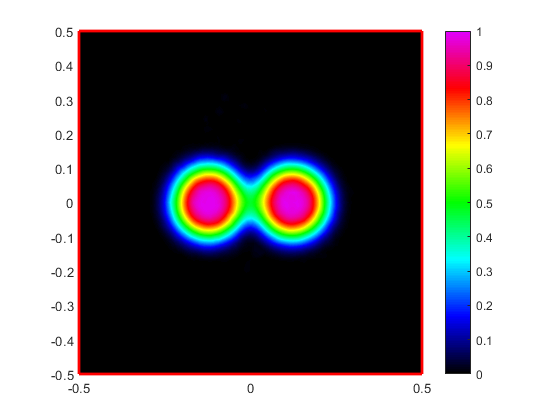}}
\subfigure{\includegraphics[width=0.24\linewidth]{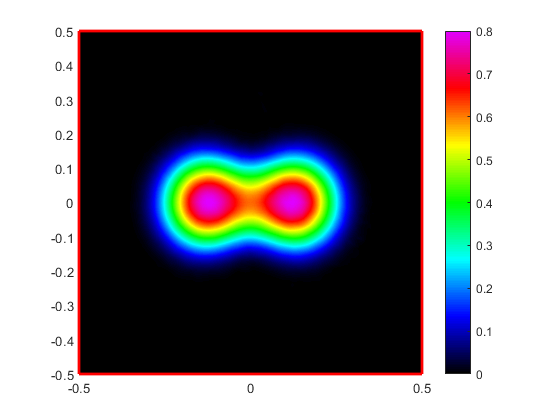}}
\subfigure{\includegraphics[width=0.24\linewidth]{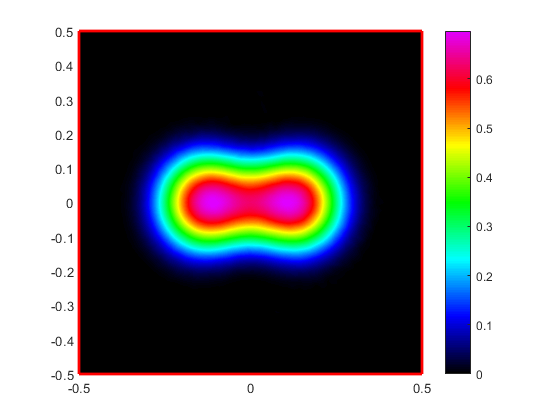}}
\subfigure{\includegraphics[width=0.24\linewidth]{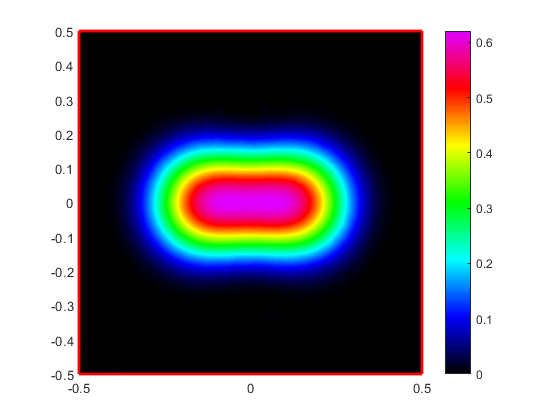}}
\subfigure{\includegraphics[width=0.24\linewidth]{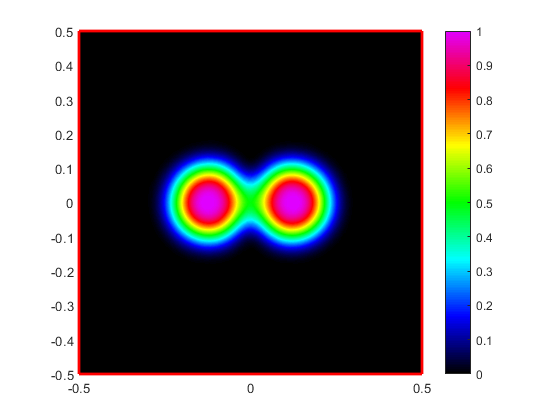}}
\subfigure{\includegraphics[width=0.24\linewidth]{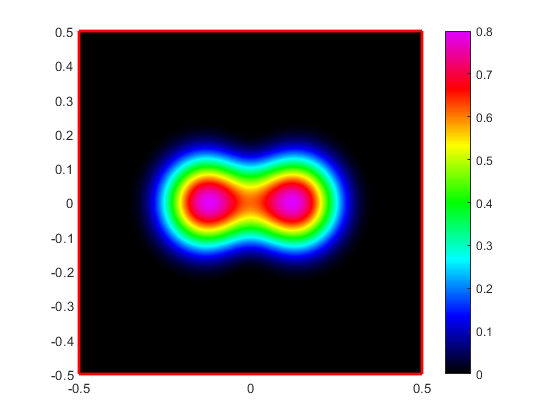}}
\subfigure{\includegraphics[width=0.24\linewidth]{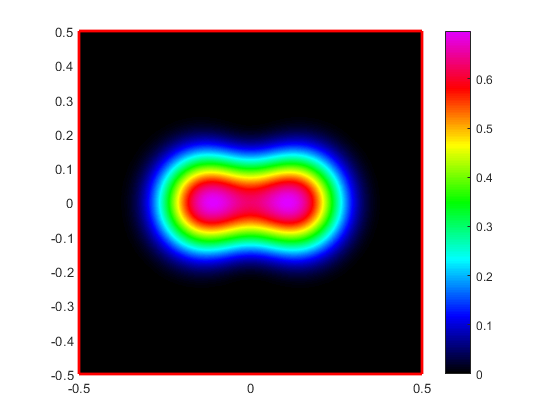}}
\subfigure{\includegraphics[width=0.24\linewidth]{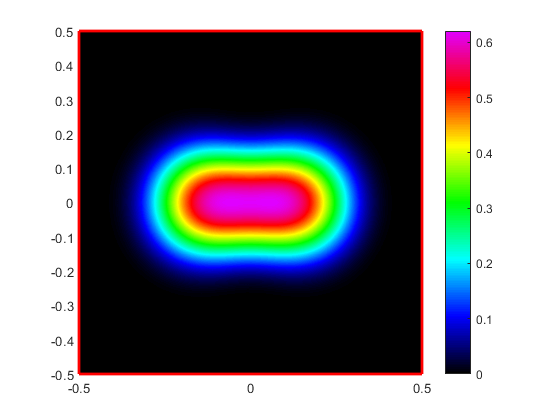}}
\subfigure{\includegraphics[width=0.24\linewidth]{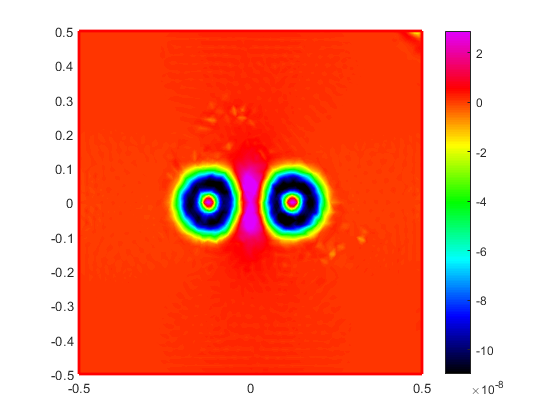}}
\subfigure{\includegraphics[width=0.24\linewidth]{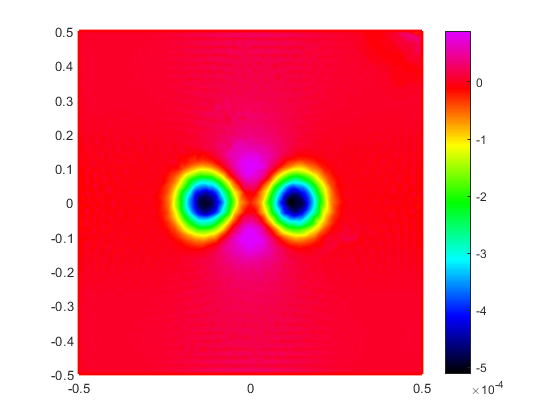}}
\subfigure{\includegraphics[width=0.24\linewidth]{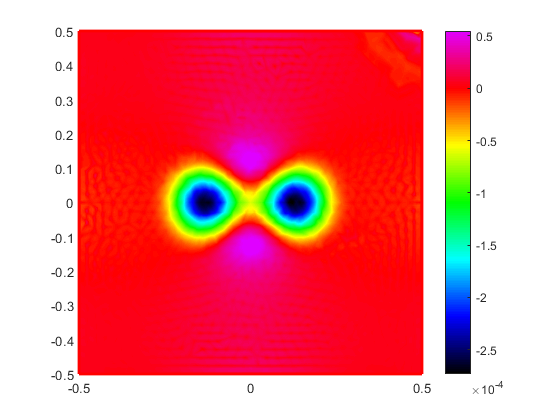}}
\subfigure{\includegraphics[width=0.24\linewidth]{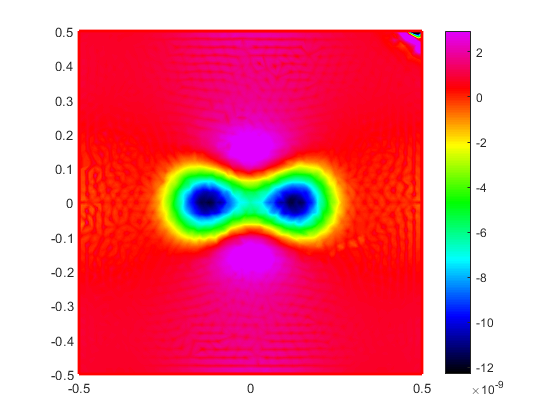}}
\subfigure{\includegraphics[width=0.24\linewidth]{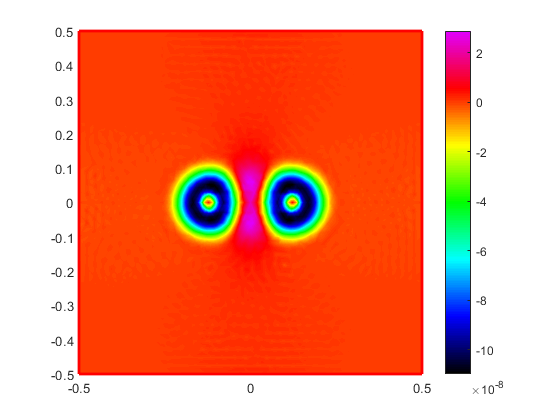}}
\subfigure{\includegraphics[width=0.24\linewidth]{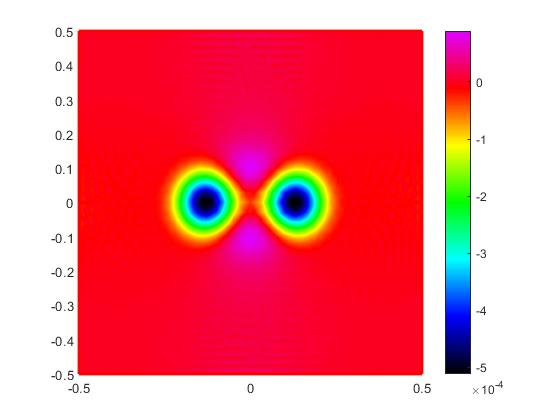}}
\subfigure{\includegraphics[width=0.24\linewidth]{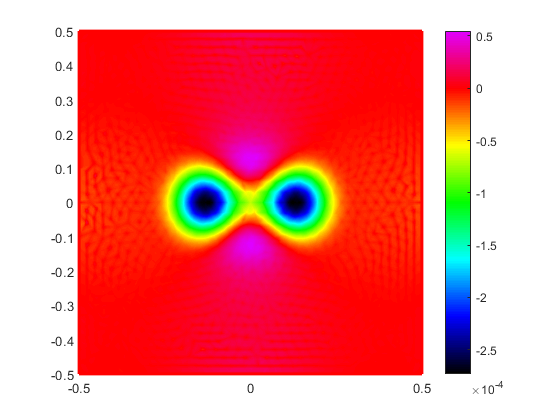}}
\subfigure{\includegraphics[width=0.24\linewidth]{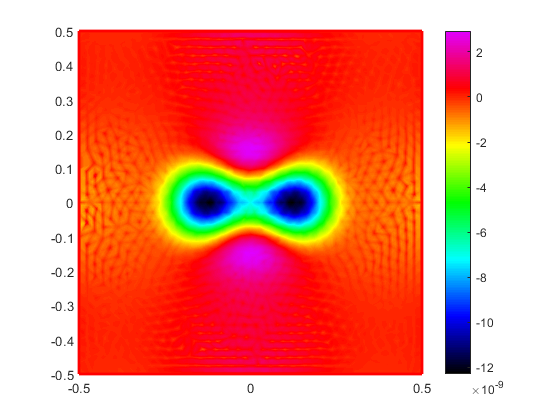}}
\caption{Predicted data by the Neural-PDE ($1^{st}$ row) and the exact data ($2^{nd}$ row) of  volume fraction $\phi$, predicted pressure $p$ ($3^{rd}$ row) and exact pressure($4^{th}$ row). The graphs of columns 1-4 represent the time states of $t_1,t_2,t_3,t_4$, respectively, where $0 \leq t_1 < t_2 < t_3 < t_4 \leq 1$.}
\label{fig:NSCH_data}
\end{figure}

This  fluid system can be derived by the energetic variational approach \citep{forster2013mathematical}. 
Here $\vu$ is the velocity and $\phi(x,y,t) \in [ 0,1 ]$ is the labeling function of the fluid phase.  $M$ is the diffusion coefficient, and $\mu$ is the chemical potential of $\phi$. Equation (\ref{eq4.4})  indicates the incompressibility of the fluid.
Solving such PDE system is notorious because of its high nonlinearity and multi-physical and coupled features. A challenge of deep learning in solving a system like this is how to process the data to improve the learning efficiency when the input matrix consists of  variables such as $\phi \in [0,1]$ with large magnitude value and variable of very small values such as $p \sim 10^{-5}$. For the Neural-PDE to better extract and learn the physical features of variables in different spatial-temporal scales, we normalized the $p$ data with a $s
igmoid$ function. We set $\delta t = 5\times 10^{-4}$. Here the training dataset is generated by the FEM solver FreeFem++~\citep{MR3043640} using a Crank-Nicolson in time $C^0$ finite element scheme. Our Neural-PDE prediction shows that the physical features of $p$ and $\phi$ have been successfully captured  with an overall MSE: $6.1631 \times 10^{-7}$ (see Figure~\ref{fig:NSCH_data}). In this example, we only coupled $p$ and $\phi$ together to show the learning ability of the Neural-PDE.  Another approach is to couple $p$, $\phi$ and the velocity $\vu$ together in the training data to predict all the related variables ($p,\phi, \vu$), which would need more normalization and regularization, techniques such as batch normalization would be helpful, please see recent research on PINN based neural network in solving such system~\citep{wight2020solving}.

\begin{figure}[t]
\centering
 \subfigure[Exact Test Dataset]{\includegraphics[width=0.49\linewidth]{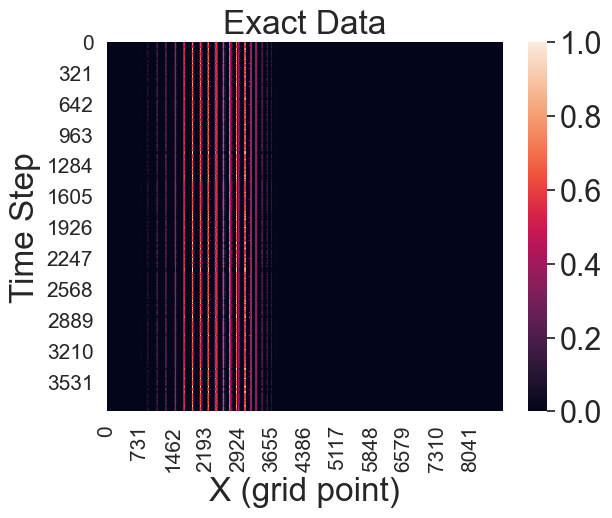}}
 \subfigure[Predicted Test Dataset]{\includegraphics[width=0.49\linewidth]{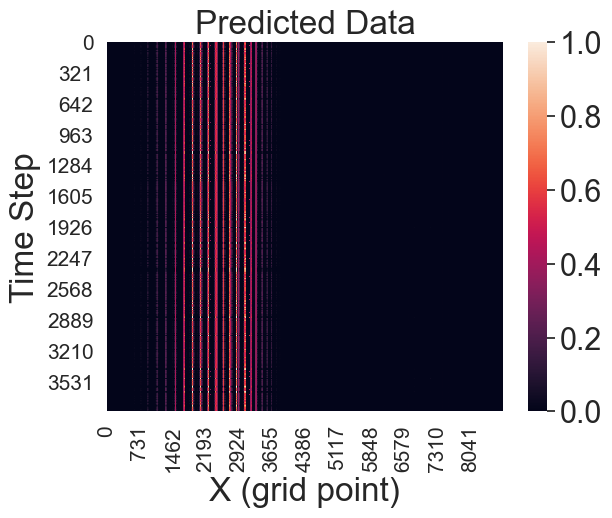}}
\subfigure[Training Metrics]{\includegraphics[width=0.45\linewidth]{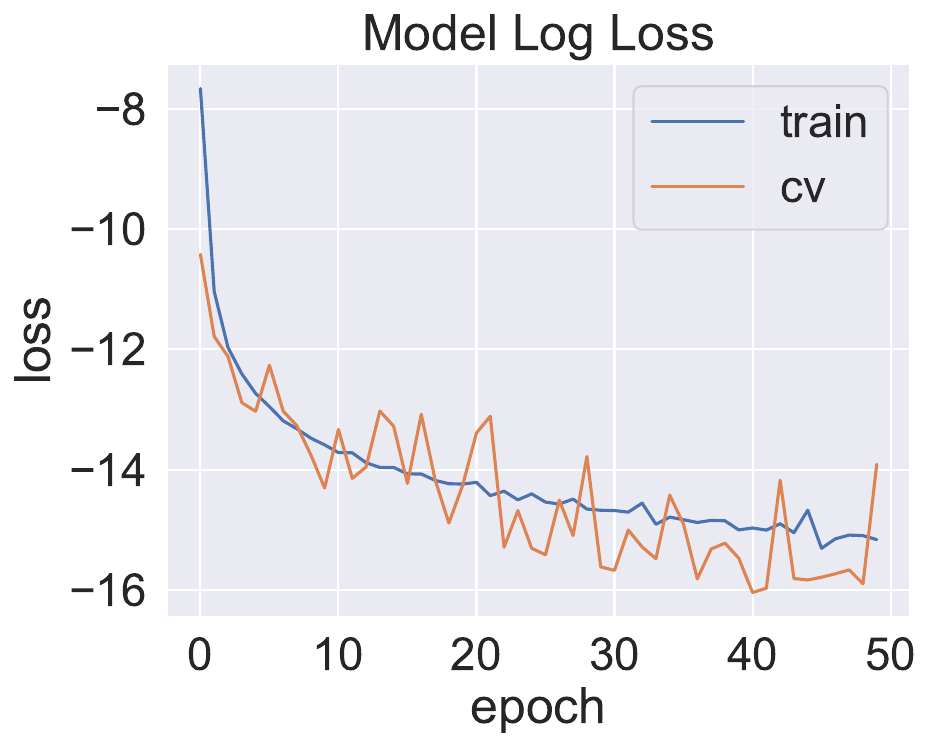}}
 \caption{Neural-PDE shows ideal prediction on Fluid System.}
 \label{fig:todo}
\end{figure}


\section{Conclusions}
In this paper, we proposed a novel sequence recurrent deep learning framework: Neural-PDE, which is capable of  intelligently filtering and learning solutions of time-dependent PDEs. One key innovation of our method  is that  the time marching method from the numerical PDEs  is applied in the deep learning framework, and the neural network is trained to explore the accurate numerical solutions for prediction.

Our experiments show that the Neural-PDE is capable of simulating from 1D to multi-dimensional  scalar PDEs  to  highly nonlinear and coupled PDE systems with their initial conditions, boundary conditions without knowing the specific forms of the equations. Solutions to the PDEs can be either continuous or discontinuous.

The state-of-the-art researches have shown the promising power of deep learning in solving high-dimensional nonlinear problems in engineering, biology and finance with efficiency in computation and accuracy in prediction. 
However, there are still  unresolved issues in applying deep learning in PDEs. For instance, the stability and convergence of the traditional numerical algorithms have been rigorously studied by applied mathematicians. Due to the high nonlinearity of the neural network system , theorems guiding stability and convergence of solutions predicted by the neural network are yet to be revealed.

Lastly, it would be helpful and interesting if one can theoretically characterize a numerical scheme from the neural network coefficients and learn the forms or mechanics from the scheme and prediction. We leave these questions for future study.
The code and data  for this paper
will become available at \url{https://github.com/YihaoHu/Neural_PDE} upon publication.

\bibliographystyle{elsarticle-num}
\bibliography{references}

\end{document}

%% file: math_commands.tex
\usepackage{amsmath,amsfonts,bm}

\def\eqref#1{equation~\ref{#1}}

\def\1{\bm{1}}

\def\rmW{{\mathbf{W}}}

\def\va{{\bm{a}}}
\def\vb{{\bm{b}}}
\def\vc{{\bm{c}}}

\def\vf{{\bm{f}}}

\def\vh{{\bm{h}}}
\def\vi{{\bm{i}}}

\def\vo{{\bm{o}}}

\def\vu{{\bm{u}}}

\def\vx{{\bm{x}}}
\def\vy{{\bm{y}}}

\def\mA{{\bm{A}}}
\def\mB{{\bm{B}}}

\DeclareMathAlphabet{\mathsfit}{\encodingdefault}{\sfdefault}{m}{sl}
\SetMathAlphabet{\mathsfit}{bold}{\encodingdefault}{\sfdefault}{bx}{n}













%% file: main.bbl
\begin{thebibliography}{10}
\expandafter\ifx\csname url\endcsname\relax
  \def\url#1{\texttt{#1}}\fi
\expandafter\ifx\csname urlprefix\endcsname\relax\def\urlprefix{URL }\fi
\expandafter\ifx\csname href\endcsname\relax
  \def\href#1#2{#2} \def\path#1{#1}\fi

\bibitem{courant1967partial}
R.~Courant, K.~Friedrichs, H.~Lewy, On the partial difference equations of
  mathematical physics, IBM journal of Research and Development 11~(2) (1967)
  215--234.

\bibitem{osher1988fronts}
S.~Osher, J.~A. Sethian, Fronts propagating with curvature-dependent speed:
  algorithms based on hamilton-jacobi formulations, Journal of computational
  physics 79~(1) (1988) 12--49.

\bibitem{leveque1992numerical}
R.~J. LeVeque, Numerical methods for conservation laws, Vol.~3, Springer, 1992.

\bibitem{cockburn2012discontinuous}
B.~Cockburn, G.~E. Karniadakis, C.-W. Shu, Discontinuous Galerkin methods:
  theory, computation and applications, Vol.~11, Springer Science \& Business
  Media, 2012.

\bibitem{thomas2013numerical}
J.~W. Thomas, Numerical partial differential equations: finite difference
  methods, Vol.~22, Springer Science \& Business Media, 2013.

\bibitem{johnson2012numerical}
C.~Johnson, Numerical solution of partial differential equations by the finite
  element method, Courier Corporation, 2012.

\bibitem{peng2020multiscale}
G.~C. Peng, M.~Alber, A.~B. Tepole, W.~R. Cannon, S.~De, S.~Dura-Bernal,
  K.~Garikipati, G.~Karniadakis, W.~W. Lytton, P.~Perdikaris, et~al.,
  Multiscale modeling meets machine learning: What can we learn?, Archives of
  Computational Methods in Engineering (2020) 1--21.

\bibitem{litjens2017survey}
G.~Litjens, T.~Kooi, B.~E. Bejnordi, A.~A.~A. Setio, F.~Ciompi, M.~Ghafoorian,
  J.~A. Van Der~Laak, B.~Van~Ginneken, C.~I. S{\'a}nchez, A survey on deep
  learning in medical image analysis, Medical image analysis 42 (2017) 60--88.

\bibitem{devlin2018bert}
J.~Devlin, M.-W. Chang, K.~Lee, K.~Toutanova, Bert: Pre-training of deep
  bidirectional transformers for language understanding, arXiv preprint
  arXiv:1810.04805 (2018).

\bibitem{lecun1998gradient}
Y.~LeCun, L.~Bottou, Y.~Bengio, P.~Haffner, Gradient-based learning applied to
  document recognition, Proceedings of the IEEE 86~(11) (1998) 2278--2324.

\bibitem{krizhevsky2012imagenet}
A.~Krizhevsky, I.~Sutskever, G.~E. Hinton, Imagenet classification with deep
  convolutional neural networks, in: Advances in neural information processing
  systems, 2012, pp. 1097--1105.

\bibitem{hinton2012deep}
G.~Hinton, L.~Deng, D.~Yu, G.~E. Dahl, A.-r. Mohamed, N.~Jaitly, A.~Senior,
  V.~Vanhoucke, P.~Nguyen, T.~N. Sainath, et~al., Deep neural networks for
  acoustic modeling in speech recognition: The shared views of four research
  groups, IEEE Signal processing magazine 29~(6) (2012) 82--97.

\bibitem{Han2018SolvingHP}
J.~Han, A.~Jentzen, E.~Weinan, Solving high-dimensional partial differential
  equations using deep learning, Proceedings of the National Academy of
  Sciences 115 (2018) 8505 -- 8510.

\bibitem{long2018pde}
Z.~Long, Y.~Lu, X.~Ma, B.~Dong, Pde-net: Learning pdes from data, in:
  International Conference on Machine Learning, 2018, pp. 3208--3216.

\bibitem{sirignano2018dgm}
J.~Sirignano, K.~Spiliopoulos, Dgm: A deep learning algorithm for solving
  partial differential equations, Journal of computational physics 375 (2018)
  1339--1364.

\bibitem{raissi2019physics}
M.~Raissi, P.~Perdikaris, G.~E. Karniadakis, Physics-informed neural networks:
  A deep learning framework for solving forward and inverse problems involving
  nonlinear partial differential equations, Journal of Computational Physics
  378 (2019) 686--707.

\bibitem{sherstinsky2020fundamentals}
A.~Sherstinsky, Fundamentals of recurrent neural network (rnn) and long
  short-term memory (lstm) network, Physica D: Nonlinear Phenomena 404 (2020)
  132306.

\bibitem{ascher1997implicit}
U.~M. Ascher, S.~J. Ruuth, R.~J. Spiteri, Implicit-explicit runge-kutta methods
  for time-dependent partial differential equations, Applied Numerical
  Mathematics 25~(2-3) (1997) 151--167.

\bibitem{chen2018neural}
R.~T. Chen, Y.~Rubanova, J.~Bettencourt, D.~K. Duvenaud, Neural ordinary
  differential equations, in: Advances in neural information processing
  systems, 2018, pp. 6571--6583.

\bibitem{huang2015bidirectional}
Z.~Huang, W.~Xu, K.~Yu, Bidirectional lstm-crf models for sequence tagging,
  arXiv preprint arXiv:1508.01991 (2015).

\bibitem{schuster1997bidirectional}
M.~Schuster, K.~K. Paliwal, Bidirectional recurrent neural networks, IEEE
  transactions on Signal Processing 45~(11) (1997) 2673--2681.

\bibitem{hochreiter1997long}
S.~Hochreiter, J.~Schmidhuber, Long short-term memory, Neural computation 9~(8)
  (1997) 1735--1780.

\bibitem{lipton2015critical}
Z.~C. Lipton, J.~Berkowitz, C.~Elkan, A critical review of recurrent neural
  networks for sequence learning, arXiv preprint arXiv:1506.00019 (2015).

\bibitem{graves2005framewise}
A.~Graves, J.~Schmidhuber, Framewise phoneme classification with bidirectional
  lstm and other neural network architectures, Neural networks 18~(5-6) (2005)
  602--610.

\bibitem{mikolov2011extensions}
T.~Mikolov, S.~Kombrink, L.~Burget, J.~{\v{C}}ernock{\`y}, S.~Khudanpur,
  Extensions of recurrent neural network language model, in: 2011 IEEE
  international conference on acoustics, speech and signal processing (ICASSP),
  IEEE, 2011, pp. 5528--5531.

\bibitem{lu2018beyond}
Y.~Lu, A.~Zhong, Q.~Li, B.~Dong, Beyond finite layer neural networks: Bridging
  deep architectures and numerical differential equations, in: International
  Conference on Machine Learning, 2018, pp. 3276--3285.

\bibitem{vaswani2017attention}
A.~Vaswani, N.~Shazeer, N.~Parmar, J.~Uszkoreit, L.~Jones, A.~N. Gomez,
  {\L}.~Kaiser, I.~Polosukhin, Attention is all you need, in: Advances in
  neural information processing systems, 2017, pp. 5998--6008.

\bibitem{forster2013mathematical}
J.~Forster, Mathematical modeling of complex fluids, Master's, University of
  Wurzburg (2013).

\bibitem{MR3043640}
F.~Hecht, \href{https://freefem.org/}{New development in freefem++}, J. Numer.
  Math. 20~(3-4) (2012) 251--265.
\newline\urlprefix\url{https://freefem.org/}

\bibitem{wight2020solving}
C.~L. Wight, J.~Zhao, Solving allen-cahn and cahn-hilliard equations using the
  adaptive physics informed neural networks, arXiv preprint arXiv:2007.04542
  (2020).

\end{thebibliography}
